\documentclass[12pt]{article}
\usepackage{eqsection,subeqnarray,indent,amsfonts,amssymb}
\usepackage{pstricks,pst-node,pst-text,pst-3d,pst-plot}
\usepackage{amsmath}     
\usepackage{bm}          
\usepackage{graphicx}
\usepackage{cite}
\usepackage{url}

\def\today{October 19, 2010 \\ revised Oct 18, 2011 \\[1mm]
                               final revision May 31, 2013}

\begin{document}

\bibliographystyle{plain}

\date{\today}

\title{Is the five-flow conjecture almost false?} 

\author{
  {\small Jesper Lykke Jacobsen${}^{1,2}$
          and Jes\'us Salas${}^{3,4}$}                   \\[1mm]
  {\small\it ${}^1$Laboratoire de Physique Th\'eorique,
                          \'Ecole Normale Sup\'erieure}  \\[-0.2cm]
  {\small\it 24 rue Lhomond, 75231 Paris, FRANCE}        \\[-0.2cm]
  {\small\tt JESPER.JACOBSEN@ENS.FR}                     \\[1mm]
  {\small\it ${}^2$Universit\'e Pierre et Marie Curie,
             4 place Jussieu, 75252 Paris, FRANCE}       \\[1mm]
  {\small\it ${}^3$Grupo de Modelizaci\'on, Simulaci\'on Num\'erica y
             Matem\'atica Industrial} \\[-2mm]
  {\small\it Universidad Carlos III de Madrid} \\[-2mm]
  {\small\it Avda.\  de la Universidad, 30}    \\[-2mm]
  {\small\it 28911 Legan\'es, SPAIN}           \\[-2mm]
  {\small\tt JSALAS@MATH.UC3M.ES}              \\[1mm]
  {\small\it ${}^4$Grupo de Teor\'{\i}as de Campos y F\'{\i}sica
             Estad\'{\i}stica}\\[-2mm]
  {\small\it Instituto Gregorio Mill\'an, Universidad Carlos III de
             Madrid}\\[-2mm]
  {\small\it Unidad Asociada al IEM-CSIC}\\[-2mm]
  {\small\it Madrid, SPAIN}           \\[4mm]
  {\protect\makebox[5in]{\quad}}  
  \\
}

\maketitle
\thispagestyle{empty}   

\begin{abstract}
  The number of nowhere zero ${\mathbb Z}_Q$ flows on a 
  graph $G$ can be shown to be a polynomial in $Q$, defining the flow
  polynomial $\Phi_G(Q)$. According to Tutte's five-flow
  conjecture, $\Phi_G(5) > 0$ for any bridgeless $G$. 
  A conjecture by Welsh that
  $\Phi_G(Q)$ has no real roots for $Q \in (4,\infty)$ was recently
  disproved by Haggard, Pearce and Royle. These authors conjectured
  the absence of roots for $Q \in [5,\infty)$.
  We study the real roots of
  $\Phi_G(Q)$ for a family of non-planar cubic graphs known as
  generalised Petersen graphs $G(m,k)$. 
  We show that the modified conjecture on real flow roots is also false, 
  by exhibiting infinitely many real flow roots
  $Q>5$ within the class $G(nk,k)$. In particular, we compute
  explicitly the flow polynomial of $G(119,7)$, 
  showing that it has real roots at $Q\approx 5.0000197675$ and 
  $Q\approx 5.1653424423$. We moreover prove that
  the graph families $G(6n,6)$ and $G(7n,7)$ possess real flow roots
  that accumulate at $Q=5$ as $n\to\infty$ (in the latter case 
  from above and below); and that 
  $Q_c(7)\approx 5.2352605291$ is an accumulation point 
  of real zeros of the flow polynomials for $G(7n,7)$ as $n\to\infty$.
\end{abstract} 


\clearpage

\newcommand{\be}{\begin{equation}}
\newcommand{\ee}{\end{equation}}
\newcommand{\ba}{\begin{subeqnarray}}
\newcommand{\ea}{\end{subeqnarray}}
\newcommand{\<}{\langle}
\renewcommand{\>}{\rangle}
\newcommand{\widebar}{\overline}
\def\reff#1{(\protect\ref{#1})}
\def\spose#1{\hbox to 0pt{#1\hss}}
\def\ltapprox{\mathrel{\spose{\lower 3pt\hbox{$\mathchar"218$}}
 \raise 2.0pt\hbox{$\mathchar"13C$}}}
\def\gtapprox{\mathrel{\spose{\lower 3pt\hbox{$\mathchar"218$}}
 \raise 2.0pt\hbox{$\mathchar"13E$}}}
\def\textprime{${}^\prime$}
\def\proof{\par\medskip\noindent{\sc Proof.\ }}
\def\qed{\hbox{\hskip 6pt\vrule width6pt height7pt depth1pt \hskip1pt}\bigskip}
\def\proofof#1{\bigskip\noindent{\sc Proof of #1.\ }}
\def\half{ {1 \over 2} }
\def\third{ {1 \over 3} }
\def\twothird{ {2 \over 3} }
\def\smfrac#1#2{\textstyle\frac{#1}{#2}}
\def\smhalf{ \smfrac{1}{2} }
%
%
\newcommand{\card}{\mathop{\rm card}\nolimits}
\renewcommand{\dim}{\mathop{\rm dim}\nolimits}
\newcommand{\Tr}{\mathop{\rm Tr}\nolimits}
\newcommand{\real}{\mathop{\rm Re}\nolimits}
\renewcommand{\Re}{\mathop{\rm Re}\nolimits}
\newcommand{\imag}{\mathop{\rm Im}\nolimits}
\renewcommand{\Im}{\mathop{\rm Im}\nolimits}
\newcommand{\sgn}{\mathop{\rm sgn}\nolimits}
\newcommand{\tr}{\mathop{\rm tr}\nolimits}
\newcommand{\diag}{\mathop{\rm diag}\nolimits}
\newcommand{\Gal}{\mathop{\rm Gal}\nolimits}
\newcommand{\mycup}{\mathop{\cup}}
\newcommand{\Arg}{\mathop{\rm Arg}\nolimits}
\def\hboxscript#1{ {\hbox{\scriptsize\em #1}} }
\def\zhat{ {\widehat{Z}} }
\def\phat{ {\widehat{P}} }
\def\qtilde{ {\widetilde{q}} }
\renewcommand{\mod}{\mathop{\rm mod}\nolimits}
\renewcommand{\emptyset}{\varnothing}

\def\scra{\mathcal{A}}
\def\scrb{\mathcal{B}}
\def\scrc{\mathcal{C}}
\def\scrd{\mathcal{D}}
\def\scrf{\mathcal{F}}
\def\scrg{\mathcal{G}}
\def\scrl{\mathcal{L}}
\def\scro{\mathcal{O}}
\def\scrp{\mathcal{P}}
\def\scrq{\mathcal{Q}}
\def\scrr{\mathcal{R}}
\def\scrs{\mathcal{S}}
\def\scrt{\mathcal{T}}
\def\scrv{\mathcal{V}}
\def\scrz{\mathcal{Z}}

\def\q{{\sf q}}

\def\Z{{\mathbb Z}}
\def\R{{\mathbb R}}
\def\C{{\mathbb C}}
\def\Q{{\mathbb Q}}
\def\N{{\mathbb N}}

\def\T{{\mathsf T}}
\def\H{{\mathsf H}}
\def\V{{\mathsf V}}
\def\D{{\mathsf D}}
\def\J{{\mathsf J}}
\def\P{{\mathsf P}}
\def\QQ{{\mathsf Q}}
\def\RR{{\mathsf R}}

\def\bsigma{{\boldsymbol{\sigma}}}
\def\bone{\bm{1}}
\def\vv{\bm{v}}
\def\uu{\bm{u}}
\def\ww{\bm{w}}

\newtheorem{theorem}{Theorem}[section]
\newtheorem{definition}[theorem]{Definition}
\newtheorem{proposition}[theorem]{Proposition}
\newtheorem{lemma}[theorem]{Lemma}
\newtheorem{corollary}[theorem]{Corollary}
\newtheorem{conjecture}[theorem]{Conjecture}
\newtheorem{result}[theorem]{Result}


\newenvironment{sarray}{
          \textfont0=\scriptfont0
          \scriptfont0=\scriptscriptfont0
          \textfont1=\scriptfont1
          \scriptfont1=\scriptscriptfont1
          \textfont2=\scriptfont2
          \scriptfont2=\scriptscriptfont2
          \textfont3=\scriptfont3
          \scriptfont3=\scriptscriptfont3
        \renewcommand{\arraystretch}{0.7}
        \begin{array}{l}}{\end{array}}

\newenvironment{scarray}{
          \textfont0=\scriptfont0
          \scriptfont0=\scriptscriptfont0
          \textfont1=\scriptfont1
          \scriptfont1=\scriptscriptfont1
          \textfont2=\scriptfont2
          \scriptfont2=\scriptscriptfont2
          \textfont3=\scriptfont3
          \scriptfont3=\scriptscriptfont3
        \renewcommand{\arraystretch}{0.7}
        \begin{array}{c}}{\end{array}}

%
%
\newcommand{\stirlingsubset}[2]{\genfrac{\{}{\}}{0pt}{}{#1}{#2}}
\newcommand{\stirlingcycle}[2]{\genfrac{[}{]}{0pt}{}{#1}{#2}}
\newcommand{\assocstirlingsubset}[3]{%
{\genfrac{\{}{\}}{0pt}{}{#1}{#2}}_{\! \ge #3}}
\newcommand{\assocstirlingcycle}[3]{{\genfrac{[}{]}{0pt}{}{#1}{#2}}_{\ge #3}}
\newcommand{\euler}[2]{\genfrac{\langle}{\rangle}{0pt}{}{#1}{#2}}
\newcommand{\eulergen}[3]{{\genfrac{\langle}{\rangle}{0pt}{}{#1}{#2}}_{\! #3}}
\newcommand{\eulersecond}[2]{\left\langle\!\! \euler{#1}{#2} \!\!\right\rangle}
\newcommand{\eulersecondgen}[3]{%
{\left\langle\!\! \euler{#1}{#2} \!\!\right\rangle}_{\! #3}}
\newcommand{\binomvert}[2]{\genfrac{\vert}{\vert}{0pt}{}{#1}{#2}}

\clearpage

%
%
\section{Introduction} \label{sec.intro} 

Given an arbitrary graph $G$ and a set of $Q$ colours, the number of 
proper vertex $Q$--colourings of $G$ is given by 
the chromatic polynomial $\chi_G(Q)$, which is indeed a polynomial in $Q$  
\cite{Birkhoff12,Whitney32a,Whitney32b}. The four-colour theorem
states that every planar graph admits a 4-vertex-colouring (i.e., 
$\chi_G(4) > 0$ for any planar graph $G$) \cite{AH77}. 

The fact that $\chi_G(Q)$ is a polynomial in $Q$, allows us to promote $Q$
from its initial definition as a positive integer to a complex 
variable $Q\in\C$. This suggests an algebraic or even analytic approach 
to the colouring problem. 
There exist many studies of the location in
$\C$ of the roots of $\chi_G(Q)$, henceforth called chromatic roots.
These studies concern either specific graphs, or all planar graphs, or
other infinite families of graphs.

Birkhoff and Lewis \cite{BL46} have made the following conjecture:
if $G$ is planar, $\chi_G(Q)>0$ for $Q \in [4,\infty)$.
Obviously, the statement of this conjecture is stronger than 
the four-colour theorem, but
unfortunately has not yet been turned into a theorem. The corresponding
result for $Q \in [5,\infty)$ has however been proved by the same 
authors \cite{BL46} (see also \cite{Woodall_97,Thomassen_97,Oxley_78}).  

Beraha and Kahane \cite{BK79} have exhibited an infinite family of
planar graphs for which $Q=4$ can be proved to be an accumulation
point of {\em complex} chromatic roots. In that sense, the four-colour
theorem is ``almost false''. Improving on this, Royle \cite{Royle05} has
proved, for a slightly different family, 
that $Q=4$ is also an accumulation point of {\em real}
chromatic roots (converging to $Q=4$ from below). 
Finally, Sokal \cite{Sokal04} has proved that for a
specific family of planar graphs (generalised $\Theta$-graphs) chromatic
roots are dense in $\C$ (except perhaps in the disc $|Q-1|<1$).

\medskip

A close cousin of the chromatic polynomial is the so-called 
{\em flow polynomial}\/ $\Phi_G(Q)$. Let $G=(V,E)$ be an arbitrary (not
necessarily planar) graph $G$ with vertex set $V$ and edge set $E$, and
let $\Gamma$ be an additive Abelian group. 
A $\Gamma$--{\em flow}\/ on $G$ is a map $\phi \colon E \to \Gamma$ 
that attributes a variable $\phi(e)$ to each edge $e \in E$, subject 
to the conservation of these variables at each vertex, with respect 
to an arbitrary chosen orientation of $E$. 
An elementary example of a flow is the current
$\phi$ in an electrical network, in which case the conservation
constraint is known as Kirchhoff's first law \cite{Kirchhoff}. 

A {\em nowhere zero $\Gamma$--flow}\/ is a $\Gamma$--flow $\phi$ such that 
$\phi(e)\neq 0$ for all $e\in E$ \cite{Tutte_book,Jaeger,Zhang}. 
If $\Gamma$ is a {\em finite}\/ Abelian group of order $Q$,
it can be shown that the number of nowhere zero $\Gamma$--flows depends only
on $Q$ (not on the {\em specific}\/ structure of the group $\Gamma$), and
it is in fact the restriction to $Q\in\N$ of a polynomial in $Q$ called
the flow polynomial $\Phi_G(Q)$ \cite{Tutte50}. 
One can then again extend the definition to $Q \in \C$
and study the location of (real or complex) flow roots.

A {\em nowhere zero $Q$--flow} of $G$ is a nowhere zero $\Z$--flow 
$\phi$ such that $|\phi(e)|\le |Q-1|$ for all $e\in E$. Tutte \cite{Tutte50}
showed that $G$ has a nowhere zero $Q$--flow if and only if it has a 
nowhere zero $\Z_Q$--flow; but these two concepts are different! 
Tutte's result immediately implies the following interesting
(but far from obvious) property of nowhere zero $\Z_Q$-flows
\cite{Tutte54}: 

\begin{proposition} \label{prop.flow.Q}
If $\Phi_G(Q)>0$ for some $Q\in\N$, then $\Phi_G(Q')>0$ for 
all integers $Q'\ge Q$. 
\end{proposition}

When $G$ is planar, one has \cite{Tutte47} the duality relation 
$\chi_{G^*}(Q) = Q \ \Phi_G(Q)$, where $G^*$ denotes the dual graph. 
In this case, the properties of $\Phi_G(Q)$ thus follow from those of 
$\chi_{G^*}(Q)$. But for non-planar $G$, the flow polynomial 
$\Phi_G(Q)$ is a genuinely new object.

It is worth stressing that the Birkhoff--Lewis theorem \cite{BL46} provides
a uniform upper bound for the real zeros of the chromatic polynomial of 
{\em all}\/  loopless planar graphs, namely $Q=5$. 
However, such an upper bound (if it actually exists at all!) is not known 
for the real zeros of the flow polynomial of arbitrary bridgeless graphs. 
(Obviously, $\Phi_G(Q)=0$ if $G$ has a bridge, because of the ``nowhere zero'' 
condition.)
The existence of such uniform upper bound and its value, if it does exist,
are long-standing open problems in Combinatorics.

Consider now arbitrary (not necessarily planar) bridgeless graphs $G$. 
Because there exist graphs not admitting a nowhere zero 4--flow,
the strongest possible results for integer and real flow roots are given,
respectively, by the following two well--known conjectures:

\begin{conjecture}%
[Tutte's five--flow conjecture \protect\cite{Tutte54,Tutte_book,Jaeger}] 
\label{conj.tutte}
For any bridgeless graph $G$, $\Phi_G(5) > 0$. 
\end{conjecture}

\noindent
{\bf Remarks.} 1. This conjecture implies that $\Phi_G(Q)>0$ for all integers
$Q\ge 5$ by Proposition~\ref{prop.flow.Q}. 

2. The Petersen graph---which is a special case $G(5,2)$ of the generalised 
Petersen graphs $G(m,k)$ to be defined in Section~\ref{sec.tm} below---has 
the flow polynomial $\Phi_{G(5,2)}(Q) = (Q-1)(Q-2)(Q-3)(Q-4)(Q^2-5Q+10)$,
which vanishes at $Q=4$. So it does not admit a nowhere zero 4--flow. 

\medskip 

\begin{conjecture}[Welsh \protect\cite{Welsh}] \label{conj.welsh}
For any bridgeless graph $G$, $\Phi_G(Q) > 0$ for $Q \in (4,\infty)$.
\end{conjecture}

It should be noted that the Welsh conjecture parallels that of
Birkhoff and Lewis for the chromatic polynomial: 
the only difference is that the endpoint $Q=4$ is included in the 
Birkhoff--Lewis conjecture for the chromatic polynomial, but not in the
Welsh conjecture for the flow polynomial.
Some results by Jackson on zero--free intervals for the
flow polynomials of cubic graphs \cite{Jackson_03,Jackson_07} also suggest 
this close parallelism between $\chi_G(Q)$ for planar $G$ and 
$\Phi_G(Q)$ for {\em arbitrary}\/ $G$. (Note that both polynomials are
evaluated at the {\em same} value of $Q$.) 

A number of weaker results have been proved over the years, notably:

\begin{theorem}[Seymour \protect\cite{Seymour81}]
\label{6_flow_thm}
For any bridgeless graph $G$, $\Phi_G(6) > 0$.
\end{theorem}

\begin{theorem}[Steinberg \protect\cite{Steinberg84}]
\label{Steinberg_thm}
For any bridgeless graph $G$ that is embeddable in the projective plane,
$\Phi_G(5)>0$. 
\end{theorem}

An immediate corollary of Seymour's theorem (using 
Proposition~\ref{prop.flow.Q}) is that $\Phi_G(Q)>0$ for all integers $Q\ge 6$.
Thus $Q=5$ is a uniform upper bound for {\em integer} flow roots. But
the above results give no clue about the existence of a uniform upper bound for
{\em real}\/ flow roots. 

\medskip

The first step into proving (or disproving) 
Conjectures~\ref{conj.tutte}--\ref{conj.welsh} consists in studying the flow
roots of ``small'' graphs.  By computing the flow roots of 
small graphs with high girth (up to $32$ vertices and girth 
at least 7),
Haggard, Pearce, and Royle \cite{HPR} have very recently found an 
explicit counterexample to the Welsh conjecture: the flow polynomial of 
the generalised Petersen graph $G(16,6)$ has two real roots larger than $Q=4$: 
$Q_1 \approx 4.0252205$, and $Q_2\approx 4.2331455$.  
However, the same authors conjectured the following modification
of Conjecture~\ref{conj.welsh}, in which 4 is replaced by 5, and the endpoint
$Q=5$ is now included in accordance with Tutte's five-flow conjecture:

\begin{conjecture}[Haggard--Pearce--Royle \protect\cite{HPR}]
\label{conj.HPR}
For any bridgeless graph $G$, $\Phi_G(Q) > 0$ for $Q \in [5,\infty)$.
\end{conjecture}

\medskip

\noindent
{\bf Remark}. 
Note that Kochol \protect\cite{Kochol06} proved that the smallest 
counterexample to Tutte's five--flow conjecture should have girth at least 9.
Note also that Jackson \cite[Corollary~39]{Jackson_03} observed, as a special
case of a more general matroidal result proved but not stated(!) by 
Oxley \cite{Oxley_78}, that if $G$ {\em and all its 3-edge-connected minors} 
have girth $\le g$, then $\Phi_G(Q)>0$ for all real $Q> g$. So any graph with
a large real flow root must either have high girth or have a 
3-edge-connected minor with high girth.

\bigskip

Even though the naive parallelism between $\chi_G(Q)$ and $\Phi_G(Q)$
has been invalidated by the above-mentioned counterexample to 
Conjecture~\ref{conj.welsh}, a related line of reasoning would be that all
these conjectures and theorems might be related by replacing $Q$ for 
the chromatic polynomial with $Q+1$ for the flow polynomial. 
Thus, the four-colour theorem \cite{AH77} ``translates'' into the 
Tutte five-flow conjecture \cite{Tutte54,Tutte_book}, and the 
Birkhoff--Lewis conjecture translates into Conjecture~\ref{conj.HPR}.
(Note that the translation of Royle's result \cite{Royle05},
showing the existence of a family of plane triangulations with real 
chromatic roots converging to 4 from below,
is consistent with the fact that Conjecture~\ref{conj.welsh} 
\cite{Welsh} is false.)

\medskip

In this paper we study the flow polynomial on the infinite family of
graphs known as the generalised Petersen graphs $G(m,k)$. 
Our main results are the following:

\begin{theorem} \label{theo.1}
The value $Q=5$ is an isolated accumulation point of real zeros of 
the flow polynomial $\Phi_G(Q)$ for the families of 
bridgeless graphs $G(6n,6)$ and $G(7n,7)$ with $n\ge 3$. Moreover:
\begin{enumerate}
\item[(a)] There is a sequence of real zeros $\{Q_n\}$ of the flow polynomials
           $\Phi_{G(6n,6)}$ that converges to $Q=5$ from below.  
\item[(b)] There is a sequence of real zeros $\{Q_n\}$ of the flow polynomials
           $\Phi_{G(7n,7)}$ that converges to $Q=5$. The sub-sequence with
           odd (resp.\/ even) $n$ converges to $Q=5$ from above (resp.\/ 
           below).
\end{enumerate}
\end{theorem}

\begin{theorem} \label{theo.2}
\mbox{}
\begin{itemize}
\item[(a)] The bridgeless graph $G(119,7)$ has flow roots at
      $Q \approx 5.00002$ and $Q \approx 5.16534$
      (where $\approx$ means ``within $10^{-5}$'').  

\item[(b)] The value $Q_c(7)\approx 5.235261$ (where $\approx$ means
      ``within $10^{-6}$'') is an 
      accumulation point of real zeros of the flow polynomials 
      $\Phi_{G(7n,7)}(Q)$. In particular, the sub-sequence for odd $n$  
      of the real zeros $\{Q_n\}$ of the flow polynomials $\Phi_{G(7n,7)}$ 
      converges to $Q_c(7)$ from below.
\end{itemize}
\end{theorem}

\bigskip

\noindent
{\bf Remark}. The largest real flow root we have explicitly found is 
$Q_0 \approx 5.1653424423$ for $G(119,7)$. 

\bigskip

Thus, the Welsh conjecture and Conjecture~\ref{conj.HPR} 
(the ``translated Birkhoff--Lewis conjecture'') are both false, 
and the Tutte five-flow conjecture is ``almost
false'' in the same sense that the four-colour theorem is  
``almost false'' \cite{BK79,Royle05}. On the other hand, 
Theorem~\ref{theo.1} includes the ``translated version'' of the existence
theorem of Royle \cite{Royle05} for real chromatic roots. 

In this work, we have considered the family of graphs $G(nk,k)$ with
$k\le 7$ and $n>2$. For each $k$, we have located the set of accumulation 
points in the complex $Q$-plane of the roots of the flow polynomial 
$\Phi_{G(nk,k)}$, as $n\to\infty$.
Most accumulation points belong to limiting curves ${\cal B}_k$; and in
particular, we are interested in locating the points $Q_c(k)$, defined 
as the largest real value where the limiting curves ${\cal B}_k$ cross 
the real axis. (These points are likely to be accumulation points as 
$n\to\infty$ of real zeros, as in Theorem~\ref{theo.2}(b); but not
always: see Section~\ref{sec.real.flow.roots}.) We have been able to 
obtain the values of $Q_c(k)$ for $k\le 11$; and the numerical 
extrapolation of these values to $k\to\infty$ yields
$\lim_{k \to\infty} Q_c(k) = Q_0 \approx 5.69$ \cite{JS_Petersen}. 
We expect that this value is the largest real accumulation point that 
one can get from the  family $G(nk,k)$.

Based on this---and on the failure of Conjectures~\ref{conj.welsh} and 
\ref{conj.HPR}---we venture the following weaker conjecture: 

\begin{conjecture} \label{conj.JS}
For any bridgeless graph $G$,
$\Phi_G(Q) > 0$ for $Q \in [6,\infty)$. 
\end{conjecture}
The disproof of Conjecture~\ref{conj.HPR} leaves basically three possibilities:

\begin{enumerate}
\item $Q \to Q+1$ {\em translation is valid}. Then Tutte's 5-flow conjecture 
      is true (because the 4-colour theorem is true) and 
      Conjecture~\ref{conj.JS} is true (because the Birkhoff--Lewis theorem 
      is true), but the  Birkhoff--Lewis conjecture is false (because 
      Conjecture~\ref{conj.HPR} is false). 

\item $Q \to Q+2$ {\em translation is valid}. Then Tutte's 5-flow conjecture 
      is false (because not every planar graph is 3--colourable), but 
      Seymour's 6-flow theorem \cite{Seymour81} corresponds to the 4-colour
      theorem. The existence of graphs with real flow roots in $(5,6)$ 
      corresponds
      to the existence of planar graphs with real chromatic roots in $(3,4)$;
      and Royle's theorem on the existence of plane triangulations with real
      chromatic roots converging to 4 from below suggests that there should
      exist graphs with real flow roots converging to 6 from below. Finally,
      the Birkhoff--Lewis conjecture and Conjecture~\ref{conj.JS} would either
      be both true or both false. 

\item {\em No translation holds: the two problems are less closely related 
      than previously thought}. In this case, Conjecture~\ref{conj.JS} 
      might well be false. Indeed, it might even be the case that there does
      not exist {\em any}\/ finite upper bound for the real flow roots of 
      general graphs; this would signal the strongest possible failure of the
      analogy between real chromatic roots of planar graphs and real flow
      roots of general graphs. 
\end{enumerate}

Note that we do in fact exhibit {\em infinitely} many flow roots in an 
interval $Q \in [4,Q_0]$ with $Q_0 >5$. This means, loosely speaking, 
that if Tutte's 5-flow conjecture is true one should look for  
a purely combinatorial proof, i.e., one that considers only integer $Q$.
This is exactly the same situation as for the 4-- and 5--colour 
theorems; they hold true even though one can find families of graphs with 
real roots approaching $Q=4$ from below \cite{Royle05}, and other families
with complex roots approaching densely to $Q=4$ \cite{BK79} 
and $Q=5$ \cite{Sokal04}.
If Conjecture~\ref{conj.JS} turns
out to be false, then it is very plausible that there exist {\em no upper
bound}\/  for real flow roots of arbitrary bridgeless graphs. 

On a more technical level, we exhibit a method for computing exactly
the flow polynomial on very large generalised Petersen graphs (which
can readily be adapted to other similar graph families). This method
relies on a transfer matrix construction similar to the one employed
in our previous work \cite{JScyclic,JStorus} on the chromatic polynomial
for graphs with periodic longitudinal boundary conditions.

The paper is organised as follows. In Section~\ref{sec.flow} we define
the flow polynomial carefully and exhibit its relation to the
$Q$-state Potts model. Building on this, we show in
Section~\ref{sec.tm} how the flow polynomial for generalised Petersen
graphs can be built by a transfer matrix construction. Our results,
given in Section~\ref{sec.res1}, are obtained by implementing this
construction on a computer and pushing the computation to as large
graphs as possible.  Note that although obtained by computational
means, the flow polynomials are exact and involve no approximation
whatsoever. In Section~\ref{sec.BKW} we introduce the Beraha--Kahane--Weiss
theorem, which plays an important role in establishing our results.
In Section~\ref{sec.res2}, we describe our analytic findings about the real
zeros of the flow polynomial for this family of graphs. 
To conclude, in Appendix~\ref{appendix.lemmas}, we prove some technical
lemmas included in the text that are essential in the proofs of the main 
results of this paper (Theorems~\ref{theo.1} and~\ref{theo.2}). In 
Appendix~\ref{sec.extra} we study some additional structural properties of
the transfer matrices. Finally, in Appendix~\ref{appendix.poly}, we give
the coefficients of the flow polynomial for the generalised 
Petersen graph $G(119,7)$. 

%
%
\section{Flow polynomial} \label{sec.flow} 
 
Let $G=(V,E)$ be a connected graph and $\Gamma$ be an Abelian group. 
Assign an arbitrary orientation to each edge $e \in E$. With respect 
to any fixed vertex $i \in V$, the edges $E_i$ incident on $i$ can then 
be characterised as either ingoing or outgoing: 
$E_i = E_i^{\rm in} \cup E_i^{\rm out}$.

A $\Gamma$--{\em flow} on $G$ is a map $\phi\colon E \to \Gamma$ that 
attributes a variable $\phi(e)$ to each edge $e \in E$, subject 
to the constraint
\be
\sum\limits_{e \in E_i^{\rm in}}  \phi(e) \;=\; 
\sum\limits_{e \in E_i^{\rm out}} \phi(e)
\label{flow_constraint}
\ee
for any $i \in V$. The edge orientation is actually immaterial in
these definitions: if one wants to change the orientation of an edge
$e_0$, it suffices to change simultaneously the sign of the flow along
that edge, $\phi(e_0) \to -\phi(e_0)$. 

A {\em nowhere zero $\Gamma$--flow}\/ is a $\Gamma$--flow $\phi$ such that
$\phi(e)\neq 0$ for all $e\in E$. If $\Gamma$ is a {\em finite}\/ Abelian 
group, we denote $\Phi_G(\Gamma)$ the number of nowhere zero 
$\Gamma$--flows on $G$.  
In particular, a {\em $\Z_Q$--flow} (resp.\ a {\em nowhere zero
$\Z_Q$--flow}) on $G$ is a map $\phi\colon E \to \{0,1,\ldots,Q-1\}$
(resp.\ $\phi\colon E \to \{1,2,\ldots,Q-1\}$) for which 
the constraint \reff{flow_constraint} is imposed modulo $Q$. 
 
Let $\Gamma$ be a finite Abelian group of order $Q$. 
Clearly, the total number of $\Gamma$-flows on $G$ is $Q^{c(E)}$, 
where for any subset $E' \subseteq E$, $c(E')$ denotes the number of 
independent cycles (cyclomatic number) in the induced graph $G'=(V,E')$.  To
obtain the number of nowhere zero $\Gamma$-flows, we first subtract for
each $e \in E$ the flows for which $\phi(e) = 0$.  Since flows with
two zero--flow edges will be subtracted off twice, these must be put
back in the sum, and proceeding by inclusion-exclusion we find \cite{Wu88}
\be
\Phi_G(\Gamma) \;=\; \sum_{E' \subseteq E} (-1)^{|E|-|E'|} Q^{c(E')} \,.
\label{NZF}
\ee
By this result, $\Phi_G(\Gamma)$ depends only on $Q$ and
is indeed the restriction to positive integers of a polynomial 
in $Q$, namely \reff{NZF}. We call \reff{NZF} the {\em flow polynomial} 
of $G$ and henceforth write it as $\Phi_G(Q)$.

Meanwhile, recall the partition function of the $Q$-state Potts
model \cite{Potts52}
\be
 Z_G(Q,v) \;=\; \sum\limits_\sigma \prod\limits_{(ij) \in E}
               e^{K \delta \big( \sigma(i),\sigma(j) \big)} \,,
\label{def_Z_Potts_spins}
\ee
where the map $\sigma \colon V \to \{0,1,\ldots,Q-1\}$ is called the spin,
and $K$ is the coupling constant. The Kronecker delta function $\delta(x,y)$
is defined by $\delta(x,y)=1$ if $x=y$, and $\delta(x,y)=0$ otherwise.
We have introduced the convenient parameter $v=e^K-1$.
By expanding the edge product and performing the sum over $\sigma$,
one recovers the partition function in the Fortuin-Kasteleyn cluster
representation \cite{FK}
\be
 Z_G(Q,v) \;=\; \sum\limits_{E' \subseteq E} v^{|E'|} Q^{k(E')} \,,
 \label{Potts}
\ee
where $k(E')$ is the number of connected components in $G'=(V,E')$.

Graph theorists will recognise in \reff{Potts} [a reparametrisation
of] the Tutte polynomial \cite{Tutte54} and interpret $\sigma$ 
in \reff{def_Z_Potts_spins} as a vertex colouring.
Proper vertex colourings, i.e., those for which adjacent vertices are
coloured differently, are obtained for $K \to -\infty$, and therefore
\be
 \chi_G(Q) \;=\; Z_G(Q,-1)
\ee
is the chromatic polynomial.

Setting instead $v=-Q$ in \reff{Potts}, and using the topological
identity
\be
  k(E') \;=\; |V| - |E'| + c(E') \,,
\ee
one establishes the connection with the flow polynomial
\be
 \Phi_G(Q) \;=\; (-1)^{|E|} Q^{-|V|} Z_G(Q,-Q) \,.
 \label{flow-Potts}
\ee
Note that $\Phi_G(Q)=0$ if $G$ contains a bridge $e_0 \in E$. Indeed,
by the constraint \reff{flow_constraint} one would have $\phi(e_0) = 0$, 
preventing the existence of a nowhere zero flow.

In the case where $G$ is planar, let $G^*$ denote the dual graph.
Recall the fundamental duality relation \cite{WuWang76} of the Potts
model partition function
\be
 Z_G(Q,v) \;=\; K\,  Z_{G^*}(Q,v^*) \,,
 \label{duality_part_func}
\ee
where $v^*$ is the dual of $v$
\be
 v \, v^* \;=\; Q \,,
 \label{duality_relation}
\ee
and the proportionality factor is
\be
 K \;=\; Q^{1-|V^*|} v^{|E|} \;=\; Q^{|V|-|E|-1} v^{|E|} \,.
 \label{duality_constant}
\ee

Noticing that $v=-Q$ is dual to $v=-1$ by \reff{duality_relation}
furnishes a relation between the flow polynomial of $G$ and the
chromatic polynomial of $G^*$. Indeed, using \reff{duality_part_func}
and \reff{duality_constant} we have \cite{Tutte47}
\be
 \chi_{G^*}(Q) \;=\; Q \ \Phi_G(Q) \,.
 \label{chromatic_flow_duality}
\ee

Alternatively, the relation \reff{chromatic_flow_duality} can be
proved by noting that there exists an obvious bijection between the
nowhere zero $\Gamma$-flows on $G$ and the proper colourings of the faces
of $G$, with the colour on one face being fixed. Indeed, let $\phi$ be
a flow on $G$. Then, turning around a vertex, each time one moves from
a face $i$ to an adjacent face $j$, if the separating edge $e$ is
seen oriented to the right (resp.\ left), its flow variable $\phi(e)$
defines the colour difference $c_j-c_i = \phi(e)$ (resp.\ $c_j-c_i =
-\phi(e)$). Starting from the face with fixed colour, these
differences define the face colouring of the whole graph. The mapping
from proper colourings to flows follows similarly.

It is useful to note that for any bridgeless 3--connected 
graph $G$, one can deduce from
\reff{Potts}/\reff{flow-Potts} that $\Phi_G(Q)$ is a polynomial in $Q$ 
of degree $|E|-|V|+1$ in $Q$, and that the first two coefficients of
$\Phi_G(Q)$ are given by
\be
\Phi_G(Q) \;=\; Q^{|E|-|V|+1} \,-\, |E| \, Q^{|E|-|V|} + \ldots
\label{eq_expansion_flow}
\ee
The first term comes from the fact that there is a unique spanning graph
$(V,E')$ in \reff{Potts}/\reff{flow-Potts} with $E'=E$. 
The second term is given by the contribution of the $|E|$ spanning subgraphs 
$(V,E')$ with $E'=E\setminus e$ for each $e\in E$, and the observation
that $k(E')=1$ since $G$ is connected and bridgeless. Notice that if we
consider the spanning graph $(V,E')$ with $E'\setminus \{e,e'\}$ for any
two distinct edges $e,e'\in E$, the 3--connectedness of $G$ guarantees that
the next term in \reff{eq_expansion_flow} will be of order $Q^{|E|-|V|-1}$, 
as in this case we also have $k(E')=1$. 

%
%
\section{Transfer matrix for flow polynomials of generalised Petersen graphs} 
\label{sec.tm} 

%
%
\subsection{Generalised Petersen graphs}

The goal of this paper is to evaluate the flow polynomial on a family
of graphs $G(m,k)$ called {\em generalised Petersen graphs} and defined as
follows: let $m,k$ be positive integers such that $m>k$.
Then $G(m,k)$ is a cubic graph with $2m$ vertices
denoted $i_p$ and $j_p$ for $p=1,2,\ldots,m$: i.e., 
\be
V(G(m,k)) \;=\; \{ i_1,\ldots,i_m,j_1,\ldots,j_m\} \,.
\ee
The edge set consists of $3m$ edges 
$(i_p j_p)$, $(i_p i_{p+1})$, $(j_p j_{p+k})$, for $p=1,2,\ldots,m$,
and with all indices considered modulo $m$: i.e., 
\be
E(G(m,k)) \;=\; \{ (i_p,j_p), (i_p i_{p+1}), (j_p j_{p+k}) \mid 1\le p\le m \}
\,.
\ee
Note that $G(m,k)$ is simple for $m\neq 2k$; but it has double edges
when $m=2k$.  
These graphs were introduced by Watkins \cite{Watkins69}.
As an example, $G(12,4)$ can be drawn as follows:
$$
%
%
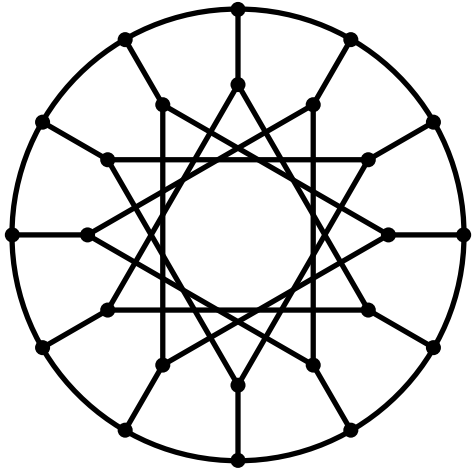
\begin{pspicture}(-3.2,-3.2)(3.2,3.2)
 \pscircle[linewidth=2pt](0.0,0.0){3.04}
 \pscircle*(3.000,0.000){0.1}
 \pscircle*(2.598,1.500){0.1}
 \pscircle*(1.500,2.598){0.1}
 \pscircle*(0.000,3.000){0.1}
 \pscircle*(-1.500,2.598){0.1}
 \pscircle*(-2.598,1.500){0.1}
 \pscircle*(-3.000,0.000){0.1}
 \pscircle*(-2.598,-1.500){0.1}
 \pscircle*(-1.500,-2.598){0.1}
 \pscircle*(0.000,-3.000){0.1}
 \pscircle*(1.500,-2.598){0.1}
 \pscircle*(2.598,-1.500){0.1}
 \pscircle*(2.000,0.000){0.1}
 \pscircle*(1.732,1.000){0.1}
 \pscircle*(1.000,1.732){0.1}
 \pscircle*(0.000,2.000){0.1}
 \pscircle*(-1.000,1.732){0.1}
 \pscircle*(-1.732,1.000){0.1}
 \pscircle*(-2.000,0.000){0.1}
 \pscircle*(-1.732,-1.000){0.1}
 \pscircle*(-1.000,-1.732){0.1}
 \pscircle*(0.000,-2.000){0.1}
 \pscircle*(1.000,-1.732){0.1}
 \pscircle*(1.732,-1.000){0.1}
 \psline[linewidth=2pt](3.000,0.000)(2.000,0.000)
 \psline[linewidth=2pt](2.598,1.500)(1.732,1.000)
 \psline[linewidth=2pt](1.500,2.598)(1.000,1.732)
 \psline[linewidth=2pt](0.000,3.000)(0.000,2.000)
 \psline[linewidth=2pt](-1.500,2.598)(-1.000,1.732)
 \psline[linewidth=2pt](-2.598,1.500)(-1.732,1.000)
 \psline[linewidth=2pt](-3.000,0.000)(-2.000,0.000)
 \psline[linewidth=2pt](-2.598,-1.500)(-1.732,-1.000)
 \psline[linewidth=2pt](-1.500,-2.598)(-1.000,-1.732)
 \psline[linewidth=2pt](0.000,-3.000)(0.000,-2.000)
 \psline[linewidth=2pt](1.500,-2.598)(1.000,-1.732)
 \psline[linewidth=2pt](2.598,-1.500)(1.732,-1.000)
 \psline[linewidth=2pt](2.000,0.000)(-1.000,1.732)
 \psline[linewidth=2pt](1.732,1.000)(-1.732,1.000)
 \psline[linewidth=2pt](1.000,1.732)(-2.000,0.000)
 \psline[linewidth=2pt](0.000,2.000)(-1.732,-1.000)
 \psline[linewidth=2pt](-1.000,1.732)(-1.000,-1.732)
 \psline[linewidth=2pt](-1.732,1.000)(0.000,-2.000)
 \psline[linewidth=2pt](-2.000,0.000)(1.000,-1.732)
 \psline[linewidth=2pt](-1.732,-1.000)(1.732,-1.000)
 \psline[linewidth=2pt](-1.000,-1.732)(2.000,0.000)
 \psline[linewidth=2pt](0.000,-2.000)(1.732,1.000)
 \psline[linewidth=2pt](1.000,-1.732)(1.000,1.732)
 \psline[linewidth=2pt](1.732,-1.000)(0.000,2.000)
\end{pspicture}
$$
The graphs $G(m,k)$ are clearly bridgeless. They are non-planar 
for all pairs $(m,k)$ except for the case $(3,2)$ and the two 
sub-families $(p,1)$ and $(2p,2)$ with $p\ge 1$. They have girth 8
for $m$ and $k$ sufficiently large. We have thus a two-parameter
family of non-planar cubic graphs with high girth, and based on
exhaustive studies of small graphs \cite{HPR} we expect this
family to produce large real flow roots.
However, it is easy to see that $Q=5$ is not a flow root:

\begin{lemma} \label{lemma.5flow}
For every generalised Petersen graph $G(m,k)$ with $m,k$ 
positive integers such that $m>k$, $\Phi_{G(m,k)}(5) > 0$. In fact, 
every graph $G(m,k)$ other than the ordinary Petersen graph
$G(5,2)$ has $\Phi_{G(m,k)}(4) > 0$. 
\end{lemma}

\proof 
It is well known \cite{Castagna72} that every generalised 
Petersen graph $G(m,k)$ (with the exception of the Petersen graph $G(5,2)$
itself), admits a Tait colouring: i.e., an edge 3--colouring such that 
at every vertex, the three incident edges take distinct colours.

It is worth noting that the definition of the generalised 
Petersen graph $G(m,k)$ in Refs.~\protect\cite{Watkins69,Castagna72} 
explicitly excludes the case $m=2k$. However, it is easy to see that 
any $G(2k,k)$ has a Tait colouring: e.g., the edges $(i_p,j_p)$ take colour $1$,
the edges $(i_p,i_{p+1})$ take alternatively colours $2$ and $3$ (as
$p$ goes from $1$ to $m$), and for each $p$, one of the double 
edges $(j_p,j_{p+k})$ takes colour $2$, and the other edge, colour $3$. 

The existence of such edge 3--colourings is equivalent, for cubic 
loopless graphs, to the existence of a nowhere zero 4--flow for the same
graph \cite[Proposition~2(b)]{Jaeger79}. Therefore, for all $G(m,k)$
except the Petersen graph $G(5,2)$, $\Phi_{G(m,k)}(4) > 0$, and 
furthermore, $\Phi_{G(m,k)}(5) > 0$ by Proposition~\ref{prop.flow.Q}.
The case $G(5,2)$ is dealt with directly: from the exact expression for
$\Phi_{G(5,2)}$ 
(see the second remark after Conjecture~\ref{conj.tutte} above), we conclude 
that $\Phi_{G(5,2)}(5)=240>0$. \qed  
 
\medskip

We shall however show that the five-flow conjecture is ``almost false'',
in the sense of Theorem \ref{theo.1}.

%
%
\subsection{Potts model transfer matrix} \label{sec.potts_tm}

We wish to evaluate $Z_{G(m,k)}(Q,v)$---of which the flow polynomial
$\Phi_{G(m,k)}(Q)$ is a special case---by a transfer matrix construction.

Contrary to an often repeated but false statement, evaluating
$Z_G(Q,v)$ by a transfer matrix construction is possible for {\em any}
graph $G$, and does not require $G$ to consist of a number of identical 
layers \cite{Bedini_10}. 
However, when $G$ does have a layered structure---as is the case
here---$Z_G(Q,v)$ can be computed by the repeated application of the
{\em same} transfer matrix.

%
%
\begin{figure}[htb]
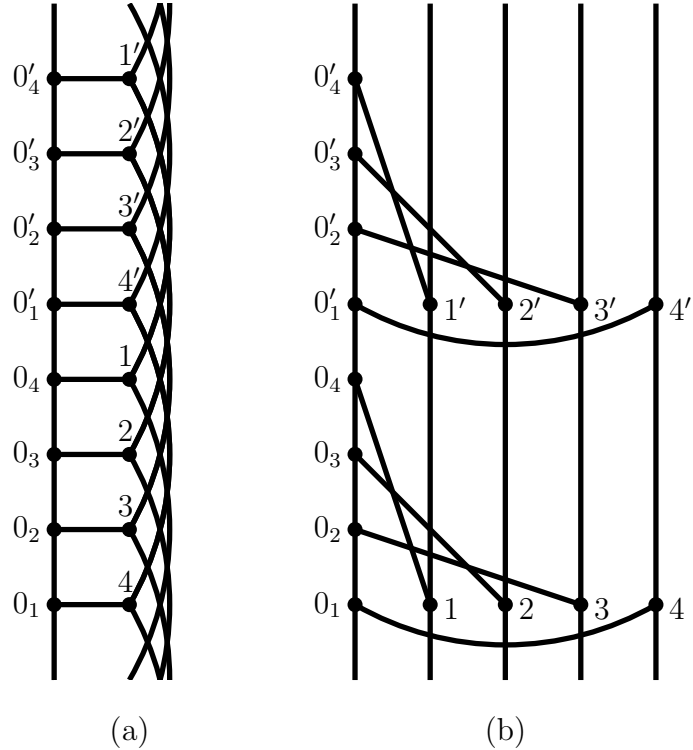

\centering
\psset{xunit=1cm}
\psset{yunit=1cm}
\pspicture*(0,-1)(10,9)
 \rput{0}(1,0){%
    \psline[linewidth=2pt](0,0)(0,9) 
    \multirput(0,1)(0,1){8}{%
       \psline[linewidth=2pt](0,0)(1,0)
       \pscircle*(0,0){0.1}
       \pscircle*(1,0){0.1}
       \psarc[linewidth=2pt](-2.4641,2){4}{-30}{30}
       \psarc[linewidth=2pt](-2.4641,-2){4}{-30}{30}
    }
    \psline*[linecolor=white](0,0)(0,-1)(2,-1)(2,0)(0,0)
    \uput[180](0,1){$0_1$}
    \uput[180](0,2){$0_2$}
    \uput[180](0,3){$0_3$}
    \uput[180](0,4){$0_4$}
    \uput[180](0,5){$0_1'$}
    \uput[180](0,6){$0_2'$}
    \uput[180](0,7){$0_3'$}
    \uput[180](0,8){$0_4'$}
    \uput[90](1,1){$4\phantom{'}$}
    \uput[90](1,2){$3\phantom{'}$}
    \uput[90](1,3){$2\phantom{'}$}
    \uput[90](1,4){$1\phantom{'}$}
    \uput[90](1,5){$4'$}
    \uput[90](1,6){$3'$}
    \uput[90](1,7){$2'$}
    \uput[90](1,8){$1'$}
    \uput[270](1,-0.3){(a)}
 }
 \rput{0}(5,0){%
   \multirput(0,0)(1,0){5}{ \psline[linewidth=2pt](0,0)(0,9) }
   \multirput(0,1)(0,4){2}{%
     \psline[linewidth=2pt](0,1)(3,0)
     \psline[linewidth=2pt](0,2)(2,0)
     \psline[linewidth=2pt](0,3)(1,0)
     \psarc[linewidth=2pt](2,3.4641){4}{240}{300}
     \multirput(0,0)(1,0){5}{\pscircle*(0,0){0.1}}
     \pscircle*(0,1){0.1}
     \pscircle*(0,2){0.1}
     \pscircle*(0,3){0.1}
    }
    \uput[180](0,1){$0_1$}
    \uput[180](0,2){$0_2$}
    \uput[180](0,3){$0_3$}
    \uput[180](0,4){$0_4$}
    \uput[180](0,5){$0_1'$}
    \uput[180](0,6){$0_2'$}
    \uput[180](0,7){$0_3'$}
    \uput[180](0,8){$0_4'$}
    \uput[350](1,1){$1$}
    \uput[350](2,1){$2$}
    \uput[350](3,1){$3$}
    \uput[350](4,1){$4$}
    \uput[350](1,5){$1'$}
    \uput[350](2,5){$2'$}
    \uput[350](3,5){$3'$}
    \uput[350](4,5){$4'$}
    \uput[270](2,-0.3){(b)}
 }
\endpspicture
 \caption{(a) Generalised Petersen graph $G(m,k)$, here with
   $k=4$.  There are $m$ layers of two vertices in the vertical
   direction, but some edges link layers at distance $k$. The boundary
   conditions are periodic in the vertical direction. (b) 
   When $m=nk$, $G(nk,k)$ can be redrawn as shown. There are
   $n=m/k$ layers of width $k+1$ vertices, each comprising a total of
   $2k$ vertices. All edges now link vertices within the same layer,
   or in two adjacent layers.}
 \label{fig:transfertrick}
\end{figure}

Let us suppose for simplicity that $m$ is a multiple of $k$: i.e., 
$m = nk$. Then the generalised Petersen graph $G(nk,k)$ can be redrawn 
as in Figure~\ref{fig:transfertrick}. This turns $G(nk,k)$ into a graph 
of $n$ identical layers of width $L=k+1$ vertices with periodic
boundary conditions in the vertical direction. (We shall henceforth
refer to this as periodic {\em longitudinal} boundary conditions, in
accordance with the fact that the transfer matrix builds up the graph
vertically.) We now claim that this implies that $Z_{G(nk,k)}(Q,v)$ can
be written as a Markov trace
\be
 Z_{G(nk,k)}(Q,v) \;=\; \Tr ({\sf T}_L)^n
 \label{Z_as_trace}
\ee
of the $n$-th power of a transfer matrix ${\sf T}_L$ to be defined 
shortly.

In general, for a layered graph of width $L$, ${\sf T}_L$ acts on
basis states $A_L$ which are set partitions of $2L$ points
$\{1,2,\ldots,L,1',2',\ldots,L'\}$. These basis states can be depicted
as {\em partition diagrams}, which are hypergraphs on $2L$ vertices,
drawn inside a rectangle with $L$ vertices (labelled
$1',2',\ldots,L'$) on top and $L$ vertices (labelled $1,2,\ldots,L$)
on bottom. Each hyperedge represents one block in the partition. A
block that contains at least one vertex from both the top and  
bottom rows is called a {\em link}. A block containing precisely one
vertex is called a {\em singleton}. The number of links in a diagram
$d$ is denoted $\ell(d)$.

The following example with $L=4$ \vspace*{2mm}
%
%
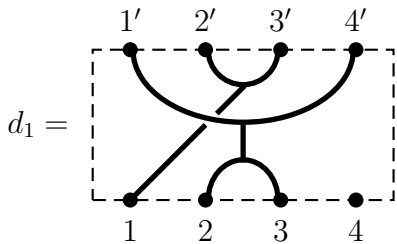
$$d_1 = 
\begin{pspicture}(-0.7,0.9)(3.7,2.5)
 \psline[linestyle=dashed](-0.5,0.0)(3.5,0.0)(3.5,2.0)(-0.5,2.0)(-0.5,0.0)
 \pscircle*(0.0,0.0){0.1}
 \pscircle*(1.0,0.0){0.1}
 \pscircle*(2.0,0.0){0.1}
 \pscircle*(3.0,0.0){0.1}
 \pscircle*(0.0,2.0){0.1}
 \pscircle*(1.0,2.0){0.1}
 \pscircle*(2.0,2.0){0.1}
 \pscircle*(3.0,2.0){0.1}
 \rput[Bc](0.0,-0.4){$1$}
 \rput[Bc](1.0,-0.4){$2$}
 \rput[Bc](2.0,-0.4){$3$}
 \rput[Bc](3.0,-0.4){$4$}
 \rput[Bc](0.0,2.4){$1'$}
 \rput[Bc](1.0,2.4){$2'$}
 \rput[Bc](2.0,2.4){$3'$}
 \rput[Bc](3.0,2.4){$4'$}
 \psellipticarc[linewidth=2pt]{-}(1.5,2.0)(1.5,1.0){180}{0}
 \psellipticarc[linewidth=2pt]{-}(1.5,0.0)(0.5,-0.5){180}{0}
 \psline[linewidth=2pt]{-}(1.5,0.5)(1.5,1.0)
 \psellipticarc[linewidth=2pt]{-}(1.5,2.0)(0.5,0.5){180}{0}
 \psline[linewidth=2pt]{-}(0.0,0.0)(1.0,1.0)
 \psline[linewidth=2pt]{-}(1.15,1.15)(1.52,1.52)
\end{pspicture} 
\vspace*{1.3cm}
$$
represents the partition $(1,2',3')(2,3,1',4')(4)$. It has two links
and one singleton.

The multiplication $d = d_2 \cdot d_1$ of two partition diagrams is
defined by stacking the diagrams vertically. Specifically, the top row
of $d_2$ becomes the top row of $d$, the bottom row of $d_1$ becomes
the bottom row of $d$, and the top row of $d_1$ is identified with the
bottom row of $d_2$. Any blocks not containing points in the top or
bottom rows of $d$ are removed in the process. This gives, for example:
\vspace*{2mm}
%
%
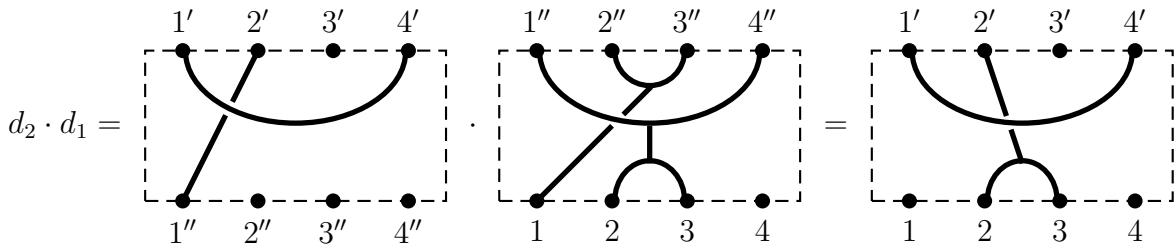
$$d_2 \cdot d_1 = 
\begin{pspicture}(-0.7,0.9)(3.7,2.5)
 \psline[linestyle=dashed](-0.5,0.0)(3.5,0.0)(3.5,2.0)(-0.5,2.0)(-0.5,0.0)
 \pscircle*(0.0,0.0){0.1}
 \pscircle*(1.0,0.0){0.1}
 \pscircle*(2.0,0.0){0.1}
 \pscircle*(3.0,0.0){0.1}
 \pscircle*(0.0,2.0){0.1}
 \pscircle*(1.0,2.0){0.1}
 \pscircle*(2.0,2.0){0.1}
 \pscircle*(3.0,2.0){0.1}
 \rput[Bc](0.0,-0.4){$1''$}
 \rput[Bc](1.0,-0.4){$2''$}
 \rput[Bc](2.0,-0.4){$3''$}
 \rput[Bc](3.0,-0.4){$4''$}
 \rput[Bc](0.0,2.4){$1'$}
 \rput[Bc](1.0,2.4){$2'$}
 \rput[Bc](2.0,2.4){$3'$}
 \rput[Bc](3.0,2.4){$4'$}
 \psellipticarc[linewidth=2pt]{-}(1.5,2.0)(1.5,1.0){180}{0}
 \psline[linewidth=2pt]{-}(0.0,0.0)(0.57,1.14)
 \psline[linewidth=2pt]{-}(0.66,1.32)(1.0,2.0)
\end{pspicture} \cdot
\begin{pspicture}(-0.7,0.9)(3.7,2.5)
 \psline[linestyle=dashed](-0.5,0.0)(3.5,0.0)(3.5,2.0)(-0.5,2.0)(-0.5,0.0)
 \pscircle*(0.0,0.0){0.1}
 \pscircle*(1.0,0.0){0.1}
 \pscircle*(2.0,0.0){0.1}
 \pscircle*(3.0,0.0){0.1}
 \pscircle*(0.0,2.0){0.1}
 \pscircle*(1.0,2.0){0.1}
 \pscircle*(2.0,2.0){0.1}
 \pscircle*(3.0,2.0){0.1}
 \rput[Bc](0.0,-0.4){$1$}
 \rput[Bc](1.0,-0.4){$2$}
 \rput[Bc](2.0,-0.4){$3$}
 \rput[Bc](3.0,-0.4){$4$}
 \rput[Bc](0.0,2.4){$1''$}
 \rput[Bc](1.0,2.4){$2''$}
 \rput[Bc](2.0,2.4){$3''$}
 \rput[Bc](3.0,2.4){$4''$}
 \psellipticarc[linewidth=2pt]{-}(1.5,2.0)(1.5,1.0){180}{0}
 \psellipticarc[linewidth=2pt]{-}(1.5,0.0)(0.5,-0.5){180}{0}
 \psline[linewidth=2pt]{-}(1.5,0.5)(1.5,1.0)
 \psellipticarc[linewidth=2pt]{-}(1.5,2.0)(0.5,0.5){180}{0}
 \psline[linewidth=2pt]{-}(0.0,0.0)(1.0,1.0)
 \psline[linewidth=2pt]{-}(1.15,1.15)(1.52,1.52)
\end{pspicture} =
\begin{pspicture}(-0.7,0.9)(3.7,2.5)
 \psline[linestyle=dashed](-0.5,0.0)(3.5,0.0)(3.5,2.0)(-0.5,2.0)(-0.5,0.0)
 \pscircle*(0.0,0.0){0.1}
 \pscircle*(1.0,0.0){0.1}
 \pscircle*(2.0,0.0){0.1}
 \pscircle*(3.0,0.0){0.1}
 \pscircle*(0.0,2.0){0.1}
 \pscircle*(1.0,2.0){0.1}
 \pscircle*(2.0,2.0){0.1}
 \pscircle*(3.0,2.0){0.1}
 \rput[Bc](0.0,-0.4){$1$}
 \rput[Bc](1.0,-0.4){$2$}
 \rput[Bc](2.0,-0.4){$3$}
 \rput[Bc](3.0,-0.4){$4$}
 \rput[Bc](0.0,2.4){$1'$}
 \rput[Bc](1.0,2.4){$2'$}
 \rput[Bc](2.0,2.4){$3'$}
 \rput[Bc](3.0,2.4){$4'$}
 \psellipticarc[linewidth=2pt]{-}(1.5,2.0)(1.5,1.0){180}{0}
 \psellipticarc[linewidth=2pt]{-}(1.5,0.0)(0.5,-0.5){180}{0}
 \psline[linewidth=2pt]{-}(1.5,0.5)(1.35,0.95)
 \psline[linewidth=2pt]{-}(1.29,1.13)(1.0,2.00)
\end{pspicture}
\vspace*{1.3cm}
$$
This diagram multiplication turns $A_L$ into an associative 
{\em partition monoid} \cite{Martin96,HR05} with identity 
$I = (1,1')(2,2')\cdots(L,L')$. Observe that
\be
 \ell(d_2 \cdot d_1) \;\le\; \min \big(\ell(d_2),\ell(d_1)\big) \,.
 \label{links_decrease}
\ee

The idea is now that these diagrams will represent the edge subset
appearing in the cluster representation \reff{Potts} of the Potts
model partition function. The factors of $v$ can be dealt with
locally, and the tricky part is to get a handle on the non-local
factors $Q$.  To this end, it is natural to associate an element of
$\C$ with each diagram, which will play the role of the Boltzmann
weight, i.e., the weight of a partially built configuration $E'$ in
\reff{Potts}. In the diagram multiplication $d = d_2 \cdot d_1$, let
$\kappa(d_1,d_2)$ 
be the number of blocks which are removed because they contain
no point in the top or bottom rows of $d$. The non-local part of
the Boltzmann weight is then $Q^{\kappa(d_1,d_2)}$.

These considerations motivate the definition of the {\em partition algebra}
\cite{Martin96,HR05} $\C A_L(Q)$ as the associative algebra over $\C$
with basis $A_L$ and multiplication defined by
\be
 d_2 \, d_1 \;=\; Q^{\kappa(d_1,d_2)} \, (d_2 \cdot d_1) \,.
\ee
The partition algebra $\C A_L(Q)$ can be represented faithfully
as an algebra of matrices in $\C^{A_L \times A_L}$ whose 
rows and columns are indexed by the partition monoid $A_L$: 
namely, the matrix $M(d)$ representing $d\in A_L$ has matrix elements
\be
M(d)_{d'' d'} \;=\; \begin{cases}
                    Q^{\kappa(d',d)} & \text{if $d'' = d \cdot d'$} \\
                    0                & \text{otherwise} 
                    \end{cases}
\ee
This is indeed the point of view that we shall take when constructing
the transfer matrix of the flow polynomial and manipulating it
explicitly (see Section~\ref{sec.algorithm}). 
The elements of $A_L$ can then be interpreted as the basis states 
of this representation. 

\medskip

\noindent
{\bf Remark.} 
With no risk of confusion, we shall therefore 
use the notation $A_L$ to refer both to the partition monoid
and to the set of basis states. However, we shall adopt a notation
that distinguishes an element ${\sf O}_L \in \C A_L(Q)$ in the
partition algebra from its corresponding matrix representation
${\cal O}_L \in \C^{A_L \times A_L}$.

\medskip

We now define a set of generators for the monoid $A_L$. These
generators will be the elementary building blocks used to define the
transfer matrix ${\sf T}_L$. Apart from the identity $I$, the
necessary generators are the {\em join operators} ${\sf J}_{ij}$ that
amalgamate the blocks containing points $i$, $j$, $i'$ and $j'$,
and the {\em detach operators} ${\sf D}_i$ that remove point $i'$ from
its block and turn it into a singleton. In the pictorial
representation this gives rise to the diagrams
\vspace*{2mm}
%
%
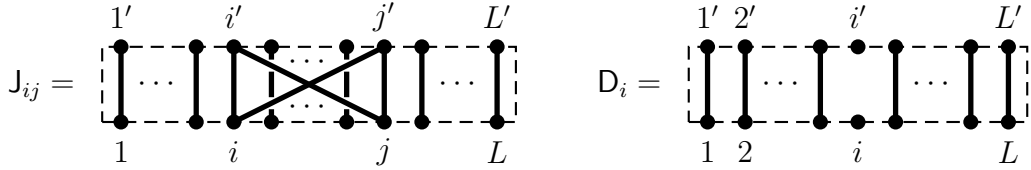
$${\sf J}_{ij} = 
\begin{pspicture}(-0.5,0.4)(5.5,1.5)
 \psline[linestyle=dashed](-0.25,0.0)(5.25,0.0)(5.25,1.0)(-0.25,1.0)(-0.25,0.0)
 \pscircle*(0.0,0.0){0.1}
 \pscircle*(1.0,0.0){0.1}
 \pscircle*(1.5,0.0){0.1}
 \pscircle*(2.0,0.0){0.1}
 \pscircle*(3.0,0.0){0.1}
 \pscircle*(3.5,0.0){0.1}
 \pscircle*(4.0,0.0){0.1}
 \pscircle*(5.0,0.0){0.1}
 \pscircle*(0.0,1.0){0.1}
 \pscircle*(1.0,1.0){0.1}
 \pscircle*(1.5,1.0){0.1}
 \pscircle*(2.0,1.0){0.1}
 \pscircle*(3.0,1.0){0.1}
 \pscircle*(3.5,1.0){0.1}
 \pscircle*(4.0,1.0){0.1}
 \pscircle*(5.0,1.0){0.1}
 \rput[Bc](0.0,-0.4){$1$}
 \rput[Bc](1.5,-0.4){$i$}
 \rput[Bc](3.5,-0.4){$j$}
 \rput[Bc](5.0,-0.4){$L$}
 \rput[Bc](0.0,1.4){$1'$}
 \rput[Bc](1.5,1.4){$i'$}
 \rput[Bc](3.5,1.4){$j'$}
 \rput[Bc](5.0,1.4){$L'$}
 \rput[Bc](0.5,0.5){$\cdots$}
 \rput[Bc](2.5,0.2){$\cdots$}
 \rput[Bc](2.5,0.8){$\cdots$}
 \rput[Bc](4.5,0.5){$\cdots$}
 \psline[linewidth=2pt]{-}(1.5,0.0)(3.5,1.0)
 \psline[linewidth=2pt]{-}(1.5,1.0)(3.5,0.0)
 \psline[linewidth=2pt]{-}(0.0,0.0)(0.0,1.0)
 \psline[linewidth=2pt]{-}(1.0,0.0)(1.0,1.0)
 \psline[linewidth=2pt]{-}(1.5,0.0)(1.5,1.0)
 \psline[linewidth=2pt]{-}(2.0,0.0)(2.0,0.17)
 \psline[linewidth=2pt]{-}(2.0,0.33)(2.0,0.67)
 \psline[linewidth=2pt]{-}(2.0,0.83)(2.0,1.0)
 \psline[linewidth=2pt]{-}(3.0,0.0)(3.0,0.17)
 \psline[linewidth=2pt]{-}(3.0,0.33)(3.0,0.67)
 \psline[linewidth=2pt]{-}(3.0,0.83)(3.0,1.0)
 \psline[linewidth=2pt]{-}(3.5,0.0)(3.5,1.0)
 \psline[linewidth=2pt]{-}(4.0,0.0)(4.0,1.0)
 \psline[linewidth=2pt]{-}(5.0,0.0)(5.0,1.0)
\end{pspicture} \qquad {\sf D}_{i} = 
\begin{pspicture}(-0.5,0.4)(4.5,1.5)
 \psline[linestyle=dashed](-0.25,0.0)(4.25,0.0)(4.25,1.0)(-0.25,1.0)(-0.25,0.0)
 \pscircle*(0.0,0.0){0.1}
 \pscircle*(0.5,0.0){0.1}
 \pscircle*(1.5,0.0){0.1}
 \pscircle*(2.0,0.0){0.1}
 \pscircle*(2.5,0.0){0.1}
 \pscircle*(3.5,0.0){0.1}
 \pscircle*(4.0,0.0){0.1}
 \pscircle*(0.0,1.0){0.1}
 \pscircle*(0.5,1.0){0.1}
 \pscircle*(1.5,1.0){0.1}
 \pscircle*(2.0,1.0){0.1}
 \pscircle*(2.5,1.0){0.1}
 \pscircle*(3.5,1.0){0.1}
 \pscircle*(4.0,1.0){0.1}
 \rput[Bc](0.0,-0.4){$1$}
 \rput[Bc](0.5,-0.4){$2$}
 \rput[Bc](2.0,-0.4){$i$}
 \rput[Bc](4.0,-0.4){$L$}
 \rput[Bc](0.0,1.4){$1'$}
 \rput[Bc](0.5,1.4){$2'$}
 \rput[Bc](2.0,1.4){$i'$}
 \rput[Bc](4.0,1.4){$L'$}
 \rput[Bc](1.0,0.5){$\cdots$}
 \rput[Bc](3.0,0.5){$\cdots$}
 \psline[linewidth=2pt]{-}(0.0,0.0)(0.0,1.0)
 \psline[linewidth=2pt]{-}(0.5,0.0)(0.5,1.0)
 \psline[linewidth=2pt]{-}(1.5,0.0)(1.5,1.0)
 \psline[linewidth=2pt]{-}(2.5,0.0)(2.5,1.0)
 \psline[linewidth=2pt]{-}(3.5,0.0)(3.5,1.0)
 \psline[linewidth=2pt]{-}(4.0,0.0)(4.0,1.0)
\end{pspicture} 
\vspace*{1.3cm}
$$
As a consequence of the above definitions, the product ${\sf D}_i d$
within $\C A_L(Q)$ produces a factor $Q$ if the point $i'$ is a
singleton in diagram $d$, and a factor $1$ otherwise. In particular,
we have 
\be
{\sf D}_i^2 \;=\; Q \, D_i \,.
\label{def_D_square}
\ee 

{}From these building blocks we can now form the operators representing
the addition of an edge to the graph that is being built up by the
transfer matrix. These are
\ba
 {\sf H}_{ij} &=& I + v {\sf J}_{ij} \,, \slabel{Hij}\\
 {\sf V}_i    &=& v I + {\sf D}_i \,. \slabel{Vi}
\label{def_H_V}
\ea
The letters {\sf H} and {\sf V} stand for horizontal and vertical, where a
horizontal edge is understood to link vertices within the same
layer of the graph (recall Figure~\ref{fig:transfertrick}), and
a vertical edge links vertices from two adjacent layers.

Inspecting Figure~\ref{fig:transfertrick}, and labelling the points as
in the figure, we can now finally define the transfer matrix (with $L=k+1$):
\begin{subeqnarray}
 {\sf T}_L &=& {\sf H}_{01}
 \left( \prod\limits_{i=k}^2 {\sf V}_0 {\sf H}_{0i} \right)
 \left( \prod\limits_{i=0}^k {\sf V}_i \right) \\[2mm]
           &=& {\sf H}_{01} {\sf V}_0 {\sf H}_{02} {\sf V}_0 {\sf H}_{03} 
               \cdots {\sf V}_0 {\sf H}_{0k} {\sf V}_k {\sf V}_{k-1}
               \cdots {\sf V}_0
 \label{TM_decomp}
\end{subeqnarray} 
Note the order of indices in the products. 

%
%
\subsection{Markov trace and eigenvalue amplitudes} \label{sec.markov}

It remains to explain the meaning of the Markov trace $\Tr$ in
\reff{Z_as_trace}. The Markov trace of any partition diagram 
$d \in A_L$ is by definition $Q^{\kappa(d)}$, where $\kappa(d)$ 
is the number of connected
components in the diagram obtained from $d$ by identifying the points
$i$ and $i'$ for all $i=1,2,\ldots,L$.  This identification
corresponds to implementing the periodic longitudinal boundary
conditions in Figure~\ref{fig:transfertrick}. The definition of the
Markov trace extends to the partition algebra $\C A_L(Q)$ by
linearity.

With this definition, \reff{Z_as_trace}, and \reff{TM_decomp},
we are in principle equipped to compute the partition function
$Z_{G(nk,k)}(Q,v)$ as a polynomial in $Q$ and $v$. A practical problem
for going to large $k$ is however that the dimension of ${\sf T}_L$,
i.e., the number of basis states $A_L$, grows very fast with $L$:
\be
 \dim \ {\sf T}_L \;=\; |A_L| \;=\; B_{2L} \,,
 \label{naive_dimension}
\ee
where $B_n$ are the Bell numbers with exponential generating function (egf)
\be
\sum\limits_{n=0}^\infty \frac{B_n z^n}{n!} \;=\; \exp({\rm e}^z - 1) \,.
\label{def_egf_bell}
\ee

Considerable progress can nevertheless be made if one takes advantage
of the structure of the partition algebra \cite{HR05}. In practical
terms this means that the number of points participating in the
partitions can be halved from $2L$ to $L$. We now explain how this
comes about.

Denote by $A_L^{(\ell)}$ the elements of the partition monoid $A_L$ with
exactly $\ell$ links, and define for $i=0,1,\ldots,L$ the set of elements 
with at most $i$ links:
\be
 {\cal A}_{L,i} \;=\; \bigcup\limits_{\ell=0}^i A_L^{(\ell)} \,.
\ee
Thanks to \reff{links_decrease}, the ${\cal A}_{L,i}$ are in fact ideals
which, moreover, constitute a filtration of the monoid:
\be
 {\cal A}_{L,0} \;\subseteq\; {\cal A}_{L,1} \;\subseteq\; \cdots
 \;\subseteq\; {\cal A}_{L,L} \;=\; A_L \,.
\ee
This implies immediately that for any element ${\sf O}_L$ (and ${\sf T}_L$ 
in particular) in the partition algebra $\C A_L(Q)$,
the corresponding matrix ${\cal O}_L \in \C^{A_L \times A_L}$ has
a block-triangular structure with respect to $\ell$. The eigenvalues
of ${\cal O}_L$ can therefore be found by restricting to
$A_L^{(\ell)}$ for $\ell = 0,1,\ldots,L$.
{}From the point of view of the matrix representation of $\C A_L(Q)$,
this restriction amounts to replacing a matrix ${\cal O}_L$ by another
matrix ${\cal O}'_L$ in which all the off-diagonal blocks have
been set to zero, i.e., ${\cal O}'_L$ is block-diagonal with 
respect to $\ell$. Although ${\cal O}'_L$ does not represent any element of
$\C A_L(Q)$, it is still a well-defined matrix and we can study
its eigenvalues, which are the same as those of ${\cal O}_L$.

In fact ${\cal O}'_L$ is block-diagonal with respect to a more refined
partition of basis states. To see this, it suffices to observe that
${\cal O}'_L$ cannot change the blocks of
the partition that contain {\em only} points from the bottom set
$\{1,2,\ldots,L\}$, since the multiplication has been defined by
acting on the top points only; nor can it amalgamate two blocks
into one, or ``abandon'' a link by failing to connect it to the top row.
Therefore, ${\cal O}'_L$ is block-diagonal,
with the blocks of ${\cal O}'_L$ being indexed by partitions
of the bottom points $\{1,2,\ldots,L\}$ together with a marking of $\ell$
of them as ``links''. Moreover, all the blocks corresponding to a given 
value of $\ell$ are {\em identical},
by virtue of the definition of the generators of $A_L$ and
the restrictions imposed when going from ${\cal O}_L$ to ${\cal O}'_L$.

As far as the determination of the eigenvalues goes, one can therefore
restrict further the basis states of $A_L^{(\ell)}$ to partitions of
the top points $\{1',2',\ldots,L'\}$ only, with precisely $\ell$
blocks (which were the links in the full partition
monoid) being marked $1,2,\ldots,\ell$ (to indicate that they are
connected, respectively, to the first, second,\ldots, $\ell$-th marked
block on the bottom row). Note that the marked blocks carry
{\em distinct} labels, since
the action of ${\cal O}'_L$ can still exchange their order (relative to
the now-forgotten fixed order of the links with respect to the bottom
points $\{1,2,\ldots,L\}$). One can then finally block-diagonalise
${\cal O}'_L$ by rearranging these restricted basis states into linear
combinations that are irreducible representations $\lambda$ of the
symmetric group $S_\ell$.

To summarise, all distinct eigenvalues of ${\cal O}_L$ can be found by
studying the irreducible representations labelled by $\ell$ and
$\lambda$. We have thus the following decomposition of the Markov trace
\be
 \Tr \ {\cal O}_L \;=\; \sum\limits_{\ell=0}^L 
                       \sum\limits_{\lambda \in S_\ell}
 \alpha_{\ell,\lambda} \ \tr_{\ell,\lambda} {\cal O}_L \,,
 \label{trace_decomp}
\ee
where now $\tr_{\ell,\lambda}$ are ordinary matrix traces. The
coefficients $\alpha_{\ell,\lambda}$ (which are polynomials in 
$Q$ as we shall see below) are eigenvalue amplitudes, which
can also be interpreted as the dimensions of the commutant of the
partition algebra.

Consider now $\lambda \in S_\ell$ through its corresponding Young diagram,
$Y(\lambda) = (\lambda_1,\lambda_2,\ldots,\lambda_\ell)$, where $\lambda_i$
is the number of boxes in the $i$--th row. If there are less than $\ell$
rows in $Y(\lambda)$, the expression is of course padded with zeros.
One then has the result \cite[Proposition 3.24]{HR05}
\be
 \alpha_{\ell,\lambda} \;=\;
   \frac{\dim \lambda}{\ell !}
   \prod\limits_{i=0}^{\ell-1} (Q-i-\lambda_{\ell-i}) \,.
 \label{eigen_amp}
\ee
We recall that the dimension $\dim \lambda$ of the representation
$\lambda$ is given by the hook formula \cite{Sagan}
\be
 \dim \lambda \;=\; \frac{\ell!}{\prod\limits_{x \in Y(\lambda)} h_x}\,,
 \label{eigen_dim}
\ee
where $h_x$ is the hook length of the box $x \in Y(\lambda)$, i.e.,
the number of boxes to its right, plus the number of boxes below it,
plus the box itself.
We shall sometimes need the total amplitude $\beta_\ell$ for a 
given number of marked blocks. This reads then
\be
 \beta_\ell \;=\; \sum\limits_{\lambda \in S_\ell} \alpha_{\ell,\lambda} \,
       \dim \lambda \,.
 \label{eigen_amp_total}
\ee

\medskip

\noindent
{\bf Remarks.}
1. A couple of other cases of decompositions of the Markov trace,
analogous to \reff{trace_decomp}--\reff{eigen_amp}, have previously
been considered in the literature. Indeed, had
the graph been planar and with periodic {\em transverse} boundary
conditions (in addition to the periodic longitudinal boundary
conditions that we assume throughout), the transfer matrix could only
have changed the cyclic order of the links, and the relevant group
would not have been $S_\ell$, but rather the cyclic group $C_\ell$. Its
representation theory leads to very different expressions
\cite[Eqs.~(1.2)/(1.3)]{Richard} 
for the analogue of \reff{trace_decomp}--\reff{eigen_amp}. A similar
remark holds for planar graphs with free transverse boundary 
conditions, in which case
the links cannot be permuted at all, and the group acting on the links
is the trivial group $E$ consisting of only the identity.
The corresponding decomposition of the Markov trace can be found in
\cite[Eqs.~(8)--(10)]{Richard_cyclic}.

2. Obviously the graph does not need to be non-planar for
\reff{trace_decomp}--\reff{eigen_amp} to be applicable.
Rather, since $E \subseteq C_\ell \subseteq S_\ell$, the two planar
cases discussed in the preceding remark can be treated in the
general non-planar formalism. But when doing so, some of the
representations $\lambda \in S_\ell$ will lead to zero eigenvalues
and/or eigenvalues corresponding to different representations
$\lambda$ will coincide. Discarding the former
representations, and summing up the amplitudes $\alpha_{\ell,\lambda}$ 
of the latter, then reproduces the results of \cite{Richard,Richard_cyclic}.

%
%
\subsection{Flow polynomial transfer matrix}

The transfer matrix that produces the flow polynomial
$\Phi_{G(nk,k)}(Q)$ can be taken simply as ${\sf T}_L$ of
Sections~\ref{sec.potts_tm}--\ref{sec.markov}, i.e., by specialising
\reff{Z_as_trace} to $v=-Q$.

One can however reduce the dimension of the relevant partition
algebra by remarking that for $v=-Q$, the vertical operator ${\sf V}_i$ 
in \reff{Vi} is a projector (up to a constant). Indeed, 
by \reff{def_D_square} one finds ${\sf V}_i^2 = (v I + {\sf D}_i)^2 =
v^2 I + (2v + Q) {\sf D_i}$, which is a multiple of ${\sf V}_i$ 
if and only if $v=-Q$. 
The normalised projector $(-Q)^{-1}{\sf V}_i$ annihilates any partition
diagram in which the point $i'$ is a singleton. Concerning the reduced
partitions of the points $\{1',2',\ldots,L'\}$ with precisely $\ell$
marked blocks---as described in Section~\ref{sec.markov}---the precise
statement is: $(-Q)^{-1} {\sf V}_i$ annihilates any reduced diagram in which
$i'$ is an un-marked singleton.

Since ${\sf V}_i$ and ${\sf V}_j$ commute for any $i,j$, the following
operator
\be
 {\sf P}_L \;=\; (-Q)^{-L} \, \prod\limits_{i=0}^{L-1} {\sf V}_i 
 \label{def_PL}
\ee
is also a projector. It annihilates any reduced diagram containing an
un-marked singleton. We can therefore replace \reff{Z_as_trace} by
\be
 Z_{G(nk,k)}(Q,-Q) \;=\; \Tr (\widetilde{\sf T}_L)^n \,,
 \label{Phi_as_trace_prev}
\ee
where $\widetilde{\sf T}_L = {\sf P}_L {\sf T}_L$ with $v=-Q$, and
consider the trace only over states without un-marked singletons.
This implies that the flow polynomial [cf.\/ \reff{flow-Potts}]
can be finally written as 
\be
 \Phi_{G(nk,k)}(Q) \;=\; (-1)^{nk} \, Q^{-2nk}\,\Tr (\widetilde{\sf T}_L)^n \,,
 \label{Phi_as_trace}
\ee
as the generalised Petersen graph $G(m,k)$ has $3m$ edges and $2m$ vertices.  
The prefactor $(-1)^{nk} \, Q^{-2nk}$ can be absorbed in the definition of the
transfer matrix: if we define
\be
\widehat{\sf T}_L \;=\; (-1)^k\, Q^{-2k} \, \widetilde{\sf T}_L \,,
\label{def_That}
\ee
then \reff{Phi_as_trace} becomes 
\be
 \Phi_{G(nk,k)}(Q) \;=\; \Tr (\widehat{\sf T}_L)^n \,.
 \label{Phi_as_trace_final}
\ee

The decomposition \reff{trace_decomp}--\reff{eigen_amp} goes through
as before, now only with the ``no un-marked singleton'' constraint
imposed on the representations labelled by $\ell$ and $\lambda$.

\medskip

\noindent
{\bf Remarks.} 
1. All the entries in the matrix $\widetilde{\sf T}_L$ are
polynomials in $Q$; but this property does not hold in general 
for the matrix elements of $\widehat{\sf T}_L$. Some of them may 
contain terms with inverse powers of $Q$.

Let us give an example for $k=3$. When we apply the transfer matrix 
$\widetilde{\sf T}_4$ to the partition $(1',2',3',4')(0,1,2,4)$, we get 
several partitions with coefficients that are polynomial in $Q$. 
In particular, we obtain the partition $(1',4')(2',3')(1,2,3,4)$ with the 
coefficient $(-Q)^5$. If we divide this polynomial by the prefactor 
$(-1)^kQ^{-2k}$ [cf.\/ \reff{def_That}], we obtain $1/Q$, which is 
{\em not}\/  a polynomial in $Q$. 

2. The structural properties of $\widetilde{\sf T}_L$ and 
   $\widehat{\sf T}_L$ are obviously the same. 

%
%
\subsection{Dimensions of representations}

Let us first consider the number of partitions $A_{m}^{(\ell)}$ 
of $m$ points with $\ell$ marked and distinguishable blocks. It is given by
\begin{equation}
 |A_m^{(\ell)}| \;=\; m! \, [z^m] \left( ({\rm e}^z-1)^\ell 
                                  \exp\left({\rm e}^z - 1\right)\right) \,,
\label{def_Akl_EGF}
\end{equation}
as is easily seen by elementary manipulations of the egf 
of the Bell numbers (the case $\ell = 0$).
Indeed, we are interested in the particular case $m=k+1$.
Using \reff{def_egf_bell}, and the fact that the Stirling numbers of the
second kind $\stirlingsubset{k}{\ell}$ (or Stirling subset numbers) 
\cite{Knuth} have the following egf \cite{Flajolet} 
\be
\stirlingsubset{k}{\ell} \;=\; \frac{k!}{\ell!} \, [z^k] 
              \left( {\rm e}^z -1 \right)^\ell \,, 
\label{def_stirling}
\ee
we can derive the following closed form for $|A_m^{(\ell)}|$:
\begin{subeqnarray}
|A_m^{(\ell)}| &=& \ell! \, \sum\limits_{p=0}^m 
\binom{m}{p} \stirlingsubset{p}{\ell} \, B_{m-p} \slabel{def_Akl1} \\
                   &=& \ell! \, \sum\limits_{s=0}^{m-\ell} 
   \binom{\ell+s}{\ell} \stirlingsubset{m}{\ell+s} \slabel{def_Akl2} 
\label{def_Akl} 
\end{subeqnarray}
where we have gone from \reff{def_Akl1} to \reff{def_Akl2} by using the
well-known expression of the Bell numbers in terms of the Stirling subset 
numbers \cite{Flajolet}
\be
B_n \;=\; \sum\limits_{s=0}^n \stirlingsubset{n}{s}\,,
\ee
and using Eq.~(6.28) of Ref.~\cite{Knuth}, valid for integers $p,n,m\ge 0$:  
\be
\stirlingsubset{n}{p+m} \, \binom{p+m}{p} \;=\; \sum\limits_k 
\stirlingsubset{k}{p} \, \stirlingsubset{n-k}{m} \, \binom{n}{k} \,. 
\ee
It is clear from \reff{def_Akl} that $\ell! \mid |A_m^{(\ell)}|$,
$|A_m^{(m)}| = m!$, and $|A_m^{(0)}| = B_m$. 

Meanwhile, the sum of the eigenvalue amplitudes for a given $\ell$
and all possible Young diagrams $\lambda$ is given by
\reff{eigen_amp}--\reff{eigen_amp_total}:
\begin{equation}
 \beta_\ell \;=\;
 \ell! \, \sum_{i=0}^\ell \frac{(-1)^i}{i!} {Q \choose \ell-i} \,. 
 \label{def.beta}
\end{equation}
This is indeed a polynomial in $Q$; it can be rewritten in terms of 
the falling factorials $Q^{\underline{i}} = \prod_{j=1}^i (Q+1-j)$ \cite{Knuth} as:
\begin{equation}
 \beta_\ell \;=\;
 \sum_{i=0}^\ell (-1)^{\ell-i} \binom{\ell}{i} \, Q^{\underline{i}} \,. 
 \label{def.beta2}
\end{equation}
The values we need in this paper are:
\begin{subeqnarray}
\beta_0 &=& 1                                                        \\
\beta_1 &=& Q   - 1                                                  \\
\beta_2 &=& Q^2 - 3Q   + 1                                           \\
\beta_3 &=& Q^3 - 6Q^2 + 8Q    -   1                                 \\
\beta_4 &=& Q^4 -10Q^3 +29Q^2  -  24Q   +   1                        \\
\beta_5 &=& Q^5 -15Q^4 +75Q^3  - 145Q^2 +  89Q   -  1                \\
\beta_6 &=& Q^6 -21Q^5 +160Q^4 - 545Q^3 + 814Q^2 -415Q    +   1      \\
\beta_7 &=& Q^7 -28Q^6 +301Q^5 -1575Q^4 +4179Q^3 -5243Q^2 +2372Q -1
\label{def.values.beta}
\end{subeqnarray}

Just as in the case of the planar partition algebra
[see Ref.~\cite{CombBLM}, in particular Eqs.~(2.16) and (2.20)]
the compatibility between the dimensions $|A_m^{(\ell)}|$
and the amplitudes $\beta_\ell$ can be expressed in the form of
a sumrule:
\begin{equation}
 \sum\limits_{\ell=0}^m \frac{1}{\ell!} \, \beta_\ell \,
 |A_m^{(\ell)}| \;=\; Q^m \,.
 \label{sumrule1}
\end{equation}
This expresses that the number of degrees of freedom per vertex of
the graph is indeed $Q$, as expected. We also have the sumrule
\be
B_{m}^{(0)} \;=\;
\sum\limits_{\ell=0}^m \frac{1}{\ell!} \, |A_m^{(\ell)}| 
\;=\; \sum\limits_{r=0}^m 2^r \, \stirlingsubset{m}{r}\,, 
\label{def_B0k}
\ee
where the integers $B_m^{(0)}$ form the sequence {\tt A001861} of
\cite{Sloane}. Their egf is
\begin{equation}
 \sum_{m=0}^\infty \frac{B_m^{(0)} z^m}{m!} \;=\;
 \exp\big( 2 ({\rm e}^z - 1) \big) \,,
\end{equation}
as can be deduced from \reff{def_Akl_EGF}.

For fixed $m$ and $\ell$, we now introduce the ``no un-marked singletons'' 
constraint. We then obtain a smaller set of partitions 
$\widetilde{A}_{m}^{(\ell)}$. We are interested in the cardinality   
of the set of partitions $\widetilde{A}_{m}^{(\ell)}$ with $m=k+1$,  
denoted by $\widetilde{N}_{k}(\ell) = |\widetilde{A}_{k+1}^{(\ell)}|$.

The number of partitions of $m$ points with $\ell$ marked and
distinguishable blocks satisfying the ``no un-marked singletons'' constraint
is 
\begin{equation}
 |\widetilde{A}_m^{(\ell)}| \;=\;
 m! \, [z^m] \left( ({\rm e}^z-1)^\ell \exp\left({\rm e}^z-1-z\right)\right) \,.
\label{def_Atilde_kl_egf}
\end{equation}
The key point is that the number of partitions with no singletons is 
given by the egf $\exp({\rm e}^z-1-z)$ \cite[p.~111]{Flajolet}. The numbers
associated with this egf are given by
\be
S_n \;=\; n! \, [z^n] \exp({\rm e}^z-1-z) \;=\; (-1)^n \sum\limits_{q=0}^n 
 \binom{n}{q} \, (-1)^{q} \, B_q \,,
\label{def_Sn}
\ee 
where the $B_n$ are the Bell numbers \reff{def_egf_bell}. Then, a closed form 
for the numbers $|\widetilde{A}_m^{(\ell)}|$ reads:
\be
|\widetilde{A}_m^{(\ell)}| \;=\; \ell! \, \sum_{p=0}^m 
\binom{m}{p}\, \stirlingsubset{p}{\ell} S_{m-p} \,.
\label{def_Atilde_kl}
\ee
This formula implies that $\ell! \mid |\widetilde{A}_m^{(\ell)}|$,
and that $|\widetilde{A}_m^{(m)}| = m!$. 

The sumrule corresponding to \reff{sumrule1} now reads
\begin{equation}
 \sum\limits_{\ell=0}^m \frac{1}{\ell!} \, \beta_\ell \,
 |\widetilde{A}_m^{(\ell)}| \;=\; (Q-1)^m \,.
\end{equation}
The fact that $Q$ has been replaced by $Q-1$ is a manifestation of the
``nowhere zero'' constraint. We also have the sumrule corresponding to 
\reff{def_B0k}
\be
\widetilde{B}_{m}^{(0)} \;=\;
\sum\limits_{\ell=0}^m \frac{1}{\ell!} \, |\widetilde{A}_m^{(\ell)}| 
\;=\; \sum\limits_{r=0}^m \binom{m}{r} B^{(0)}_r (-1)^{m-r}\,, 
\ee
where the integers $B_m^{(0)}$ are given by \reff{def_B0k}. 
Their egf is
\begin{equation}
 \sum\limits_{m=0}^\infty \frac{\widetilde{B}_m^{(0)} z^m}{m!} \;=\;
 \exp\left( 2 ({\rm e}^z - 1) -z \right) \,.
\end{equation}

Notice that the difference 
$|A_m^{(\ell)}| - |\widetilde{A}_m^{(\ell)}|$ gives the 
number of partitions of the set $\{1,2,\ldots,m\}$ with $\ell$ marked
points and with at least one un-marked singleton. These partitions do
not contribute to the final result, as they are associated to null
eigenvalues.  

\medskip

\noindent
{\bf Remark}. For simplicity, as $L=k+1$, we will consider hereafter 
the bottom-row (resp.\/ top-row) points labelled as $\{0,1,\ldots,k\}$ 
(resp.\/ $\{0',1',\ldots,k'\}$). The monoid $A_L$ will contain the partitions 
of the set $\{0,1,\ldots,k,0',1',\ldots,k'\}$.  

%
%
\section{Flow polynomial for the generalised Petersen graphs} 
   \label{sec.res1} 

\subsection{General theory} \label{sec.res1.gen}

Let us start with the simplest case $k=1$. The graph $G(n,1)$ is isomorphic 
to a cyclic ladder of width $2$. From the known Potts-model partition function
\cite[and references therein]{Tutte_sq}, one can easily derive 
\be
\Phi_{G(n,1)}(Q) \;=\; (Q^2-3Q+1) (-1)^n +(Q-1)(Q-3)^n + (Q-2)^n \,,
\ee
where the eigenvalues $\mu=-1,Q-3,Q-2$ correspond to the sectors with
$\ell=2,1,0$ links, respectively.  

Let us now focus on $k\ge 2$. 
The flow polynomial \reff{Phi_as_trace_final} for the generalised 
Petersen graph $G(nk,k)$ can be written using \reff{trace_decomp} as  
\be
\Phi_{G(nk,k)}(Q) \;=\; \sum\limits_{\ell=0}^{k+1} 
   \sum\limits_{\lambda\in S_\ell} \alpha_{\ell,\lambda} 
   \tr_{\ell,\lambda} (\widehat{\sf T}_{k+1})^n \,,  
\label{def_Phi_G_as_trace_ok}
\ee
where the amplitudes $\alpha_{\ell,\lambda}$ are given by \reff{eigen_amp}. 
This formula is the most general one. In terms of the non-zero eigenvalues 
$\mu_{k,\ell,\lambda,s}$ of the transfer matrix $\widehat{\sf T}_{k+1}$, 
it reads:
\be
\Phi_{G(nk,k)}(Q) \;=\; \sum\limits_{\ell=0}^{k+1}
   \sum\limits_{\lambda\in S_\ell} \alpha_{\ell,\lambda}
   \sum\limits_{s=1}^{\widetilde{N}_k(\ell,\lambda)} 
   \mu_{k,\ell,\lambda,s}^n \,, 
\label{def_Phi_G}
\ee
where $\widetilde{N}_k(\ell,\lambda)$ is given by 
\be
\widetilde{N}_k(\ell,\lambda) \;=\; \widetilde{N}_k(\ell) \, 
                                    \frac{\dim \lambda}{\ell!} 
\label{def_dim_T_L_l_lambda0}
\ee
(see \reff{def_dim_T_L_l_lambda} in the proof of Lemma~\ref{lemma.dim.ntK}), 
and $\widetilde{N}_k(\ell) = |\widetilde{A}_{k+1}^{(\ell)}|$ 
[cf.\/ \reff{def_Atilde_kl_egf}/\reff{def_Atilde_kl}]. 

The flow polynomial $\Phi_{G(nk,k)}$ is obtained in 
\reff{def_Phi_G_as_trace_ok} as a linear combination of ordinary matrix traces 
with definite coefficients $\alpha_{\ell,\lambda}$ given by \reff{eigen_amp}.
This is {\em all}\/  that we need to compute $\Phi_{G(nk,k)}$ rigorously 
from the various diagonal blocks of $\widehat{\sf T}_{L}$.
It is worth stressing that Eq.~\reff{def_Phi_G} holds true irrespective 
of whether some eigenvalues happen to be identical or not. 

To simplify the notation, we will denote by $\widehat{\sf T}_{k+1,\ell}$ the
diagonal block of the full transfer matrix $\widehat{\sf T}_{k+1}$ 
corresponding to partitions with exactly $\ell$ links. We will denote by 
$\widehat{\sf T}_{k+1,\ell,\lambda}$ the diagonal block of 
$\widehat{\sf T}_{k+1,\ell}$ corresponding to the irreducible representation
$\lambda$ of the group $S_\ell$. Similar notation will be used for the  
diagonal blocks of the transfer matrix $\widetilde{\sf T}_{k+1}$. Then,  
Eq.~\reff{def_Phi_G_as_trace_ok} can be rewritten as:
\be
\Phi_{G(nk,k)}(Q) \;=\; \sum\limits_{\ell=0}^{k+1}
   \sum\limits_{\lambda\in S_\ell} \alpha_{\ell,\lambda}
   \tr (\widehat{\sf T}_{k+1,\ell,\lambda})^n \,,
\label{def_Phi_G_as_trace_ok2}
\ee
where $\tr$ is an ordinary trace. The dimension of the matrix 
$\widehat{\sf T}_{k+1,\ell,\lambda}$ is given by \reff{def_dim_T_L_l_lambda0}.

We have symbolically computed {\em all} blocks 
$\widehat{\sf T}_{k+1,\ell,\lambda}$ for $1\le k\le 7$, $0\le \ell\le k+1$,
and all representations $\lambda\in S_\ell$. Therefore, we can compute 
the {\em exact} flow polynomial $\Phi_{G(nk,k)}$ for $1\le k\le 7$ by using 
Eq.~\reff{def_Phi_G_as_trace_ok2}, as all the elements involved are exactly
known.  

The dimension of the blocks for $k=6,7$ is in some cases 
very large, and the computation of the traces in \reff{def_Phi_G_as_trace_ok2}
is very memory-- and CPU-consuming even for modest values of $n$. 
For instance, the block for $k=7$, $\ell=3$, and $\lambda=(2,1)$ has 
dimension $14\,364$.

Fortunately, the blocks $\widehat{\sf T}_{k+1,\ell,\lambda}$ have for any
$k\ge 1$, any $1\le \ell \le k+1$, and any $\lambda\in S_\ell$
an additional internal structure. This is given by the following lemma 
proved in Appendix~\ref{appendix.lemmas}: 

\begin{lemma} \label{lemma.new1}
Fix $k\ge 1$. Then for any $0\le \ell\le k+1$, and any irreducible
representation $\lambda \in S_\ell$, the diagonal block 
$\widehat{\sf T}_{k+1,\ell,\lambda}$ can be written as an 
upper-block-triangular matrix when the basis vectors are ordered appropriately:
\be
\widehat{\sf T}_{k+1,\ell,\lambda} \;=\; \left( 
 \begin{array}{cc}
   \widehat{\sf D}_{k+1,\ell,\lambda}  & \widehat{\sf S}_{k+1,\ell,\lambda}\\ 
   {\sf 0}  & \widehat{\sf T}_{k+1,\ell,\lambda}^{(\text{nt})} 
  \end{array} \right) \,,
  \label{eq.lemma.new1}
\ee
where $\widehat{\sf D}_{k+1,\ell,\lambda}$ is a diagonal matrix with all its
diagonal elements equal to $\mu_{k,k+1}=(-1)^k$, whose rows and 
columns are indexed by partitions of $\{0,1,\ldots,k\}$ with no un-marked 
singletons, $\ell$ marked clusters, and vertex 0 is a marked singleton.
\end{lemma} 

\medskip

\noindent
{\bf Remarks.} 
1. The eigenvalue $\mu_{k,k+1}=(-1)^k$ will be called ``trivial'' in 
the following. The non-trivial eigenvalues come from the 
blocks $\widehat{\sf T}_{k+1,\ell,\lambda}^{(\text{nt})}$ [hence the 
superscript ``(nt)''].  

2. It follows from the description of $\widehat{\sf D}_{k+1,\ell,\lambda}$ 
that for $\ell=0$ the block $\widehat{\sf D}_{k+1,0}$ has zero dimension 
(i.e., all eigenvalues coming from the sector $\ell=0$ are non-trivial), 
and that for $\ell=k+1$ the non-trivial block 
$\widehat{\sf T}_{k+1,k+1,\lambda}^{(\text{nt})}$ has zero dimension 
(i.e., all eigenvalues coming from the sector $\ell=k+1$ are trivial).
 
3. The traces $\tr (\widehat{\sf T}_{k+1,\ell,\lambda})^n$ are trivially 
written as:
\be
\tr (\widehat{\sf T}_{k+1,\ell,\lambda})^n \;=\; 
 (-1)^{kn} \,  \dim \widehat{\sf D}_{k+1,\ell,\lambda}  \;+\; 
\tr (\widehat{\sf T}_{k+1,\ell,\lambda}^{(\text{nt})})^n \,.
\label{def_reduced_traces}
\ee
In this way we significantly reduce the burden of the computation. 
For instance,  the non-trivial block for $k=7$, $\ell=3$, and $\lambda=(2,1)$ 
has dimension $11\, 816$, compared to $14\, 364$ for the whole matrix 
$\widehat{\sf T}_{k+1,\ell,\lambda}$. 

4. There is of course a similar upper-block-triangular decomposition for the 
matrix $\widetilde{\sf T}_{k+1,\ell,\lambda}$. The trivial eigenvalues are 
given in this case by $(-1)^k \, (-1)^k Q^{2k} = Q^{2k}$ 
[cf., \reff{def_That}]. 

5. In Appendix \ref{sec.extra}, we shall prove an explicit formula for 
$\dim \widehat{\sf D}_{k+1,\ell,\lambda}$ for arbitrary $k,\ell,\lambda$ 
(Lemma~\ref{lemma.mult.trivi}). But we stress that this result plays no role
in our proof of Theorems~\ref{theo.1} and~\ref{theo.2}, since for $1\le k\le 7$ 
we have determined these dimensions by explicit computation.     

\medskip

As described in Section~\ref{sec.algorithm} below,
we have {\em exactly} computed {\em all} non-trivial blocks 
$\widehat{\sf T}_{k+1,\ell,\lambda}^{(\text{nt})}$, as well as the 
dimensions $\dim \widehat{\sf D}_{k+1,\ell,\lambda}$ for $1\le k\le 7$,
$0\le \ell\le k+1$, and $\lambda \in S_\ell$. Then, we have computed 
the flow polynomials $\Phi_{G(nk,k)}$ for $1\le k\le 7$ and selected values of
$n\ge 1$ by using \reff{def_Phi_G_as_trace_ok2}/\reff{def_reduced_traces}.  
In particular, Eq.~\reff{def_Phi_G_as_trace_ok2} and Lemma~\ref{lemma.new1}
are the {\em essential} elements in our method to compute the flow
polynomials on large graphs $G(nk,k)$, establishing in particular
Theorem~\ref{theo.2}(a). The practical implementation of this method is 
explained in detail in the next section.  

Now we will focus on those results that we need to prove Theorems~\ref{theo.1}
and~\ref{theo.2}(b). Both theorems are based on the Beraha--Kahane--Weiss 
theorem (see Section~\ref{sec.BKW}), which concerns analytic 
functions of the form \reff{def_fn0}. In particular, we should show that 
the eigenvalues  $\mu_{k+1,\ell,\lambda,s}$ [cf., \reff{def_Phi_G}] 
coming  from the transfer matrices $\widehat{\sf T}_{k+1,\ell,\lambda}$
satisfy all the hypotheses of this theorem. First we need to prove that 
the amplitudes $\alpha_{\ell,\lambda}$ [cf., \reff{eigen_amp}] and the 
eigenvalues $\mu_{k+1,\ell,\lambda,s}$ are analytic functions of $Q$ in 
some domain $D$ of the complex $Q$-plane. Indeed, the amplitudes are 
polynomials in $Q$, hence analytic functions of $Q$ in the whole complex plane. 
The eigenvalues $\mu_{k+1,\ell,\lambda,s}$ are algebraic functions of $Q$; 
they are thus analytic in the whole complex $Q$-plane, except at the 
branch cuts. So we can choose $D$ to be any connected open set of $\C$ not
containing a branch cut.  

In addition, there is a ``no-degenerate-dominance'' condition 
requiring that there must not exist two eigenvalues 
$\mu_i$ and $\mu_j$ with $i\neq j$ [the labels $i,j$ are just a 
shorthand for our indices $k+1,\ell,\lambda,s$] such that 
1) $\mu_i \equiv e^{i\delta} \mu_j$ for some real constant $\delta$, and 
2) the region $D_i \subseteq D$ where these two eigenvalues dominate 
(i.e., $|\mu_i| = |\mu_j| \ge |\mu_k|$ for all $k$) has a nonempty 
interior. The first step to check that this condition holds is to 
find out if there are equal eigenvalues among the $\mu_{k+1,\ell,\lambda,s}$. 
If two or more eigenvalues $\mu_{i_1},\ldots,\mu_{i_p}$ are {\em exactly} equal 
{\em for all} $Q$, then we can amalgamate them into a single term $\mu_i$ 
of the sum \reff{def_fn0}, and absorb the multiplicity into the corresponding 
amplitude $\alpha_i$. So each family of {\em equal} eigenvalues can be 
considered to be a single eigenvalue for the purpose of checking the 
``no-degenerate-dominance'' condition.

We have already proved that the trivial eigenvalue $\mu_{k,k+1}=(-1)^k$ 
appears in all blocks $\widehat{T}_{k+1,\ell,\lambda}$ with $\ell\ge 1$. 
Thus, we may include this eigenvalue only once in the sum \reff{def_fn0} 
arising in the Beraha--Kahane--Weiss theorem, and absorb the multiplicity in 
the amplitude $\alpha_{k,k+1}$. We now consider the non-trivial blocks 
$\widehat{T}_{k+1,\ell,\lambda}^{(\rm nt)}$. The exact symbolic computation 
of these blocks for $1\le k\le 7$ reveals that there is an exact degeneracy 
of eigenvalues when $\ell=k$. This is the content of the following lemma
proved in Appendix~\ref{appendix.lemmas}:
    
\begin{lemma} \label{lemma.new2}
Fix $k\ge 1$ and $\ell=k$. Then, for all irreducible representations 
$\lambda \in S_\ell$, there are $k$ eigenvalues $\mu_{k,k,s}$ in the 
non-trivial diagonal block $\widehat{\sf T}_{k+1,k,\lambda}^{(\rm nt)}$. 
Each of these eigenvalues $\mu_{k,k,s}$ has multiplicity $\dim\lambda$.
\end{lemma}
Thus, we may include each eigenvalue $\mu_{k,k,s}$ only once in the sum 
\reff{def_fn0}, and absorb the multiplicity in the amplitude $\alpha_{k,k,s}$.

Finally, we have to check that, for each fixed value of 
$k \in \{1,2,\ldots,7\}$, the eigenvalues $\mu_{k,\ell,\lambda,s}$ 
(for $0\le\ell\le k-1$), $\mu_{k,k,s}$ (for $\ell=k$) and 
$\mu_{k,k+1}=(-1)^k$ satisfy the 
``no-degenerate-dominance'' condition. This can be achieved by numerically 
computing the values of {\em all} these eigenvalues at a generic value of $Q$, 
and finding that there are no two eigenvalues with the same absolute value. 
In our case, we choose $Q = \pi + i \sqrt{3}$. Therefore we conclude that: 

\begin{lemma} \label{lemma.new4}
Fix $k\in\{1,2,\ldots, 7\}$. Then the non-trivial eigenvalues 
$\mu_{k,\ell,\lambda,s}$ (for $0\le\ell\le k-1$) and $\mu_{k,k,s}$ 
(for $\ell=k$), and the trivial eigenvalue $\mu_{k,k+1}=(-1)^k$ satisfy:
\begin{enumerate}  
\item For every $Q$ except perhaps a finite set, the eigenvalues are all
      distinct.
\item The ``no-degenerate-dominance'' condition holds.
\end{enumerate}
\end{lemma}

\medskip

\noindent
{\bf Remark.} We do not know how to prove the extension of this result to
$k\ge 8$, since our proof for $1\le k\le 7$ is by explicit 
computation. We nevertheless conjecture that Lemma~\ref{lemma.new4} holds 
true for all $k\ge 1$.

\medskip

Lemmas~\ref{lemma.new2} and~\ref{lemma.new4} are essential 
for proving that our eigenvalues satisfy the hypotheses of the 
Beraha--Kahane--Weiss theorem. This theorem is the starting point 
for proving Theorems~\ref{theo.1} and~\ref{theo.2}(b).

%
%
\subsection{Practical procedure} \label{sec.algorithm}

We have written a {\sc perl} script to compute the {\em symbolic}\/
transfer matrix ${\sf T}_{k+1}$ using ideas similar to those already 
explained in \cite{JScyclic,JStorus}. For $1\le k\le 4$, we
have checked our programs using {\sc Mathematica}. Further checks were
performed with code written in {\sc C} that allows us to numerically
compute the leading eigenvalue for given values of $k,\ell$, and 
$\lambda = (\ell)$ = the completely symmetric irreducible 
representation of $S_\ell$. 

The first step is to obtain the relevant diagonal blocks 
$\widetilde{\sf T}_{k+1,\ell,\lambda}$ of the transfer matrix 
$\widetilde{\sf T}_{k+1}$. 
We first fix the value of $\ell$ ($0\le \ell \le k+1$)
and a bottom-row configuration compatible with the chosen value of $\ell$ 
and the ``no un-marked singletons'' condition.

\medskip

\noindent 
{\bf Remark.} Our choice for the bottom-row partition is the simplest one. 
We take the partition $\{\{0,1,\ldots,k\}\}$ for $\ell=0$, 
the partition $\{\{\overline{0},\overline{1},\ldots,\overline{k}\}\}$ 
for $\ell=1$, the partition 
$\{\{\overline{0},\overline{1},\ldots,\overline{k-1}\},\{\overline{k}\}\}$ 
for $\ell=2$,\ldots, and the partition 
$\{ \{\overline{0}\},\{\overline{1}\}, \ldots, \{\overline{k}\}\}$ for 
$\ell=k+1$. The overline over a site means that this site (and the block
it belongs to) is connected to a block of the top-row partition by a link. 

\medskip

We then determine the basis of the relevant partition space of dimensionality 
$\widetilde{N}_k(\ell)$ [cf.~\reff{def_Atilde_kl}]. 
Indeed, the result does not depend on the chosen bottom-row 
partition.

\medskip

\noindent
{\bf Remark.} 
This statement is true if we explicitly mark $\ell$ blocks of the 
bottom-row partition, and leave un marked the rest of the blocks (if any). 
If we instead chose not to do this marking, we would arrive at a basis of 
dimensionality  $p\widetilde{N}_k(\ell)$ for some integer $p\ge 2$.   
This would happen e.g., for $k=3$ and $\ell=1$, if we chose the 
bottom-row partition $\{\{0,1\},\{2,3\}\}$ without explicitly saying which 
block is marked. Therefore, for a given top-row connectivity, e.g.,  
$\{\{0',1',2'\},\{\overline{3'}\}\}$, there would correspond {\em two} 
distinct partitions of the full set $\{0,\ldots,3,0',\ldots,3'\}$: namely,
$(0',1',2')(3',0,1)(2,3)$ and $(0',1',2')(0,1)(3',2,3)$. Each 
eigenvalue of the transfer matrix would then be repeated $p=2$ times. 
This extra factor comes obviously from 
the two ways we can mark one of the two blocks in $\{\{0,1\},\{2,3\}\}$.  
We stress that with the choice of the preceding remark  
this problem can never occur.

\medskip

We now choose an irreducible representation $\lambda$ of the symmetric group
$S_\ell$ of dimensionality $\dim\lambda$. To obtain the relevant diagonal
block $\widetilde{\sf T}_{k+1,\ell,\lambda}$ corresponding to $\lambda$, 
we simply take as our basis vectors those linear combinations of the 
``standard'' basis vectors with the appropriate properties under $S_\ell$. 
Finally, in order to ``extract'' the trivial eigenvalues $Q^{2k}$ for 
$\ell\ge 1$, we exploit the structure of 
$\widetilde{\sf T}_{k+1,\ell,\lambda}$ given by Lemma~\ref{lemma.new1}.
We thus obtain the non-trivial block 
$\widetilde{\sf T}_{k+1,\ell,\lambda}^{(\rm nt)}$ and the dimension 
$\dim \widetilde{\sf D}_{k+1,\ell,\lambda}$.
We then compute the powers $(\widetilde{T}^{(nt)}_{k+1,\ell,\lambda})^n$
and their traces, for the desired values of $n$. 

For small values of $k$, this computation can be  performed
symbolically for any not-too-large value of $n$,
using a symbolic algebraic manipulator program, such as {\sc Mathematica}. 
However, for larger values of $k$ (say, $k=6,7$) this is not feasible, as 
we have blocks of dimension as large as $11\,816$ (for $k=7$, $\ell=3$, 
and $\lambda=(2,1)$), and the symbolic computation of the powers of such 
large blocks is extremely time-- and memory--consuming, beyond our current 
computer capabilities.

A key issue in the subsequent analysis is the dependence of the 
traces $\tr (\widehat{\sf T}_{k+1,\ell,\lambda}^{(\rm nt)})^n$ on $Q$, $n$, 
and $k$. The needed information is given by the following lemma (the proof 
can be found in Appendix~\ref{appendix.lemmas}):

\begin{lemma} \label{lemma.degree}
Let $k\ge 1$, $0\le \ell \le k+1$, and $\lambda$ be an irreducible 
representation of $S_\ell$. Then, for each $n\ge 0$,   
$\tr (\widehat{\sf T}_{k+1,\ell,\lambda}^{(\rm nt)})^n$ 
[cf.\/ \protect\reff{eq.lemma.new1}] is a polynomial in $Q$ of degree
at most $n[k+\min(1-\ell,0)]$. 
\end{lemma}

\noindent
Our computation then made use of the following tricks: 
\begin{itemize}
 \item By Lemma~\ref{lemma.degree}, we know that the traces 
       $\tr (\widehat{\sf T}_{k+1,\ell,\lambda}^{(\rm nt)})^n$
       are  polynomials in $Q$ of degree at most
       $d = n[k + {\rm min}(1-\ell,0)]$. 
       Therefore, it suffices to compute the {\em evaluation}\/ of  
       each trace at $d+1$ integer values of $Q\neq 0$,
       and then reconstruct the corresponding polynomial using 
       Lagrange's interpolation method. In order to check the result, 
       we always compute at least $d+2$ values of the trace.  
       Please note that we compute the evaluation of the trace
       $\tr (\widehat{\sf T}_{k+1,\ell,\lambda}^{(\rm nt)})^n$
       by first computing the evaluation of the 
       trace $\tr (\widetilde{\sf T}_{k+1,\ell,\lambda}^{(\rm nt)})^n$, and then
       multiplying the result by the factor $(-1)^{nk}Q^{-2nk}$ 
       [cf.\/ \reff{def_That}]. This is why we take $Q\neq 0$.

 \item {\em Not}\/ all the entries of the transfer matrix 
       $\widetilde{\sf T}_{k+1,\ell,\lambda}^{(nt)}$ 
       are integers for integer values
       of $Q$; rather they are rational numbers.
       As we want to perform the trace computation with 
       (infinite-precision) integer arithmetic, for each value of $Q$
       we multiplied the matrix by the minimum (positive) 
       integer value such that all entries are integers. After the 
       computation is done, we reconstructed the true solution by 
       dividing by the appropriate factor. 

 \item The integers involved in the actual calculations are very large. 
       Therefore, we compute the value of the 
       trace for a given value of $Q\neq 0$ using modular arithmetic 
       for a given set of prime numbers $p \le 65\,521 \ltapprox 2^{16}$ 
       (we need up to 65 different primes). 
       We then reconstruct the value using the Chinese remainder theorem
       using infinite-precision arithmetic in {\sc Mathematica}. We 
       always use at least one more prime than needed, in order to
       check the result. To accelerate the computation of the trace
       using modular arithmetic, we use a program written in {\sc C}.  

 \item For $1\le k\le 6$, we are able to compute the traces for many different
       values of $n$. However, for $k=7$, the computation is so demanding,
       that we have focused on powers of the type $n=2^q +1$ with integer 
       $q\ge 1$. The reason why we consider odd powers for $k=7$ will become 
       clear in Section~\ref{sec.res2}. 
\end{itemize} 

Once the traces are computed, we can form the flow polynomial using 
\reff{def_Phi_G_as_trace_ok2}.
Notice that for $1 \le k \le 7$, everything in this
formula is exactly known. The zeros of the flow polynomials are then obtained 
using the program {\sc MPSolve} \cite{MPSolve,Bini00}. This software 
has the advantage that if one requests the zeros with 50--digit precision
(as in our case), the results are guaranteed to have at least such precision. 

\bigskip

\noindent
{\bf Remark}. To give a clear idea of what has been achieved, consider
the case $k=7$, and more specifically the computation of $\Phi_G$ for
the graph $G = G(17k,k)$, which is the largest computation undertaken
in this work. Applying naively \reff{TM_decomp} within the
diagrammatic basis would imply computing the 17th power of 
${\sf T}_{k+1}$, which according to (\ref{naive_dimension}) is a matrix of
dimension $B_{2(k+1)} = 10\,480\,142\,147$ whose entries are
polynomials in $Q$ of degree at most $2k=14$. The decomposition of 
${\sf T}_{k+1}$ and use of the ``no un-marked singleton'' constraint has
reduced the computation to the sum over $31$ blocks, the largest of
which has dimension $11\,816$. Even with these tricks, the computation
took around six months calendar time, using 50--80 processors, 
corresponding to some 30 years of CPU time.

%
%
\subsection{Additional checks}

Because our results are derived using software, we have performed some
tests in order to ensure that the results are correct. First of all, 
for the smallest members of each family $G(nk,k)$, we have computed the
flow polynomial using three different software programs: {\sc Maple},
the program {\sc Tutte} developed by Haggard, Pierce and Royle \cite{HPR},
and the program {\sc Tutte} developed by Bedini and Jacobsen 
\cite{Bedini_10}. 
The pairs $(k,n)$ for which the checks have been performed are shown in 
Table~\ref{table_checks}. In all cases, the agreement with our 
transfer-matrix computations is perfect.

%
%
\def\kk{\phantom{1}}
\begin{table}[t]
\centering
\begin{tabular}{|r|c|c|c|}
\hline\hline
$k$ & {\sc Maple} & HPR {\sc Tutte} \cite{HPR} &
                    BJ  {\sc Tutte} \cite{Bedini_10} \\
\hline
1  & $1\le n \le 15$   & $1\le n \le 50$   & $1\le n \le 50$   \\
2  & $1\le n \le 7\kk$ & $1\le n \le 25$   & $1\le n \le 25$   \\
3  & $1\le n \le 5\kk$ & $1\le n \le 11$   & $1\le n \le 12$   \\
4  & $1\le n \le 4\kk$ & $1\le n \le 8\kk$ & $1\le n \le 7\kk$ \\
5  & $1\le n \le 3\kk$ & $1\le n \le 6\kk$ & $1\le n \le 5\kk$ \\
6  & $1\le n \le 2\kk$ & $1\le n \le 5\kk$ & $1\le n \le 4\kk$ \\
7  & $1\le n \le 2\kk$ & $1\le n \le 4\kk$ & $1\le n \le 4\kk$ \\
\hline\hline
\end{tabular}
\caption{ \label{table_checks}
Tests performed on our transfer-matrix computations of the flow polynomial
for the generalised Petersen graphs $G(kn,k)$. For each value of $k$ in the
interval $1\le k \le 7$, we show the values of $n$ for which we have
computed $\Phi_{G(kn,k)}$ using a) {\sc Maple} (second column),
b) the {\sc Tutte} code developed by Haggard, Pierce, and Royle
\cite{HPR} (third column), and c) the {\sc Tutte} code
developed by Bedini and Jacobsen \cite{Bedini_10} (fourth column).
In all cases, the agreement between these computations and our transfer-matrix
results is perfect.
}
\end{table}
%
%

For the cubic graphs $G(m,k)$ that we are considering, we may
improve on \reff{eq_expansion_flow} by adding a few more terms:
\begin{eqnarray}
 \Phi_{G(m,k)}(Q) &=& Q^{|E|-|V|+1} - |E| Q^{|E|-|V|} +
 \left( \frac{|E|(|E|-1)}{2} - |V| \right) Q^{|E|-|V|-1} \nonumber \\
 &-& \left( \frac{|E|(|E|-1)(|E|-2)}{6} - |V|(|E|-2) \right) Q^{|E|-|V|-2}
 + \ldots \,.
\end{eqnarray}
In this expression, the coefficient of $Q^{|E|-|V|-1}$ arises from
two contributions in which $E \setminus E'$ is respectively two edges,
and three edges all incident on the same vertex. The contributions to
the coefficient of $Q^{|E|-|V|-2}$ are slightly more complicated to
characterise. Inserting $|V|=2m$ and $|E|=3m$ we obtain
\be
 \Phi_{G(m,k)}(Q) \;=\; Q^{m+1} \left( 1 - \frac{3m}{Q} + 
   \frac{m}{2} \, \frac{9m-7}{Q^2}
 - \frac{m}{2}\, \frac{(3m-2)(3m-5)}{Q^3} + \ldots \right) \,.
 \label{perturbation_improved}
\ee
We have checked that for all the graphs $G(nk,k)$ we have considered,
the flow polynomials obtained from the procedure outlined above
indeed satisfy \reff{perturbation_improved}.

There are some theorems that give us some information about the location of
the real zeros of the flow polynomial. The first theorem applies to a 
general bridgeless graph, while the second one applies only to cubic graphs  
(i.e., it is valid for $G(nk,k)$ with $n>2$): 

\begin{theorem}[Wakelin \protect\cite{Wakelin}; see also %
  Refs.~\protect\cite{Edwards_98,Jackson_09}]
Let $G$ be a bridgeless graph with $|V|$ vertices, $|E|$ edges, $b$ blocks,
and no isolated vertices. Then:
\begin{itemize}
 \item $\Phi_G(Q)$ is non-zero with sign $(-1)^{|E|-|V|+1}$ for 
       $Q\in(-\infty,1)$.
 \item $\Phi_G(Q)$ has a zero of multiplicity $b$ at $Q=1$.
 \item $\Phi_G(Q)$ is non-zero with sign $(-1)^{|E|-|V|+b+1}$ for 
       $Q\in(1,\smfrac{32}{27}]$.
\end{itemize}
\end{theorem}

\begin{theorem}[Jackson \protect\cite{Jackson_03,Jackson_07}]
Let $G$ be a 3--connected cubic graph with $|V|$ vertices and $|E|$ edges. 
Then:
\begin{itemize}
 \item $\Phi_G(Q)$ is non-zero with sign $(-1)^{|E|-|V|}$ for $Q\in(1,2)$.
 \item $\Phi_G(Q)$ has a zero of multiplicity $1$ at $Q=2$.
 \item $\Phi_G(Q)$ is non-zero with sign $(-1)^{|E|-|V|+1}$ for
       $Q\in(2,\delta)$, where $\delta\approx 2.546$ is the flow root
       of the cube in the interval $(2,3)$ [i.e., the zero in this 
       interval of $Q^3 -9Q^2 +29Q-32$]. 
\end{itemize}
\end{theorem}

\bigskip

The generalised Petersen graphs $G(nk,k)$ satisfy $|V|=2nk$, $|E|=3nk$, 
and $b=1$. Therefore, the above theorems imply that, for any $n> 2$, 
$\Phi_G(Q)$ has only two simple real roots in the interval $(-\infty,\delta)$,  
namely $Q=1$ and $Q=2$. The sign for $Q\in(-\infty,1)\cup (2,\delta)$ is that 
of $(-1)^{nk+1}$, and it has the opposite sign for $Q\in(1,2)$. For large
enough $Q>0$, the sign of $\Phi_G(Q)$ is always positive. Therefore,
for even $k$ and any $n\ge 3$, or for odd $k$ and even $n\ge 4$, 
this implies the existence of a real zero in $[\delta,\infty)$.

\medskip

\noindent
{\bf Remark.} The lower bound is sharp, as the cube is isomorphic to the 
generalised Petersen graph $G(4,3)$, which has a zero at $Q=\delta$. 

\medskip

We have explicitly checked that all the computed flow polynomials 
$\Phi_{G(nk,k)}$  have only two simple roots in the whole
interval $(-\infty,\delta)$, namely $Q=1,2$. Furthermore, for all even 
(resp.\/ odd) $k$, and all $n\ge 3$ (resp.\/ all even $n\ge 4$), the 
polynomial $\Phi_{G(nk,k)}$ has at least one root in $[\delta,\infty)$.

If $G$ is a cubic graph, one can easily see whether $Q=3$ is a 
flow root or not \cite[Proposition~6.4.2]{Diestel}: 

\begin{theorem}\label{theo.diestel}
A bridgeless cubic graph $G$ has a nowhere zero 3--flow if and only 
if it is bipartite.
\end{theorem}

\noindent
As the graphs $G(nk,k)$ are bipartite if and only if $k$ is odd and $n$ is
even, Theorem~\ref{theo.diestel} implies that $\Phi_{G(nk,k)}(Q)$ has
at least one factor $Q-3$ whenever $G(nk,k)$ is not bipartite: namely, when
$k$ is even or $n$ is odd (or both). When $G(nk,k)$ is bipartite,
then $\Phi_{G(nk,k)}(3)> 0$. We have checked these facts
in all the flow polynomials we have explicitly computed in this work.
 
%
%
\section{The Beraha--Kahane--Weiss theorem} \label{sec.BKW}

A central role in the subsequent analysis is played by a theorem
on analytic functions due to 
Beraha, Kahane and Weiss (BKW) \cite{BKW75,BKW78,BK79,BKW80}
and generalised slightly by Sokal \cite{Sokal04}.
The situation is as follows:
let $D$ be a domain (connected open set) in the complex plane,
and let $\alpha_1,\ldots,\alpha_M,\mu_1,\ldots,\mu_M$ ($M \ge 2$)
be analytic functions on $D$, none of which is identically zero.
For each integer $n \ge 0$, define
\be
   f_n(z)   \;=\;   \sum\limits_{k=1}^M \alpha_k(z) \, \mu_k(z)^n
   \;.
\label{def_fn0}
\ee
We are interested in the zero sets
\be
   \scrz(f_n)   \;=\;   \{z \in D \colon\;  f_n(z) = 0 \}
\ee
and in particular in their limit sets as $n\to\infty$:
\begin{eqnarray}
   \liminf \scrz(f_n)   & = &  \{z \in D \colon\;
   \hbox{every neighbourhood $U \ni z$ has a nonempty intersection} \nonumber\\
      & & \qquad \hbox{with all but finitely many of the sets } \scrz(f_n) \}
   \\[4mm]
   \limsup \scrz(f_n)   & = &  \{z \in D \colon\;
   \hbox{every neighbourhood $U \ni z$ has a nonempty intersection} \nonumber\\
      & & \qquad \hbox{with infinitely many of the sets } \scrz(f_n) \}
\end{eqnarray}
Let us call an index $k$ {\em dominant at $z$}\/ if
$|\mu_k(z)| \ge |\mu_l(z)|$ for all $l$ ($1 \le l \le M$);
and let us write
\be
   D_k  \;=\;  \{ z \in D \colon\;  k \hbox{ is dominant at } z  \}
   \;.
\ee
Then the limiting zero sets can be completely characterised as follows:

\begin{theorem}%
[Beraha--Kahane--Weiss \protect\cite{BKW75,BKW78,BK79,BKW80,Sokal04}]
   \label{theo.BKW}
Let $D$ be a domain in $\C$,
and let $\alpha_1$, $\ldots$, $\alpha_M$, $\mu_1,\ldots,\mu_M$ ($M \ge 2$)
be analytic functions on $D$, none of which is identically zero.
Let us further assume a ``no-degenerate-dominance'' condition:
there do not exist indices $k \neq k'$
such that $\mu_k \equiv \omega \mu_{k'}$ for some constant $\omega$
with $|\omega| = 1$ and such that $D_k$ ($= D_{k'}$)
has nonempty interior.
For each integer $n \ge 0$, define $f_n$ by
$$
   f_n(z)   \;=\;   \sum\limits_{k=1}^M \alpha_k(z) \, \mu_k(z)^n
   \;.
\label{def_fn}
$$
Then $\liminf \scrz(f_n) = \limsup \scrz(f_n)$,
and a point $z$ lies in this set if and only if either:
\begin{itemize}
   \item[(a)]  There is a unique dominant index $k$ at $z$,
       and $\alpha_k(z) =0$;  or
   \item[(b)]  There are two or more dominant indices at $z$.
\end{itemize}
\end{theorem}
Note that case (a) consists of isolated points in $D$,
while case (b) consists of curves
(plus possibly isolated points where
all the $\mu_k$ vanish simultaneously).
Henceforth we shall denote by $\scrb$ the locus of points
satisfying condition (b).

We shall often refer to the functions $\mu_k$ as ``eigenvalues'',
and to the $\alpha_k$ as ``amplitudes'',
because that is exactly how they arise in the transfer matrix formalism.

In Ref.~\cite[p.~55]{BKW80}, Beraha, Kahane, and Weiss give (without proof)
the following corollary, concerning the convergence of {\em real}
roots of $f_n$ to {\em real isolated} limiting points,
based on their proof of Theorem~\ref{theo.BKW}:

\begin{corollary}\label{coro.BKW}
Assume the hypotheses of Theorem~\ref{theo.BKW}. Let $z_0$ be a real 
isolated limiting point, and suppose that the functions $f_n$, the dominant 
eigenvalue $\mu_\star$ and its coefficient $\alpha_\star$ are all real in an
interval $(z_0-\epsilon,z_0 + \epsilon)$ for some $\epsilon > 0$ [with 
of course $\alpha_\star(z_0)=0$] and suppose further than 
$\alpha_\star'(z_0)\neq 0$.
Then $z$ is the limit of a real sequence $\{z_n\}$, defined for all 
sufficiently large $n$, for which $f_n(z_n)=0$. 
\end{corollary} 

The proof is simple: for any sufficiently small $\epsilon > 0$, 
$\alpha_\star(z_0\pm\epsilon)$ are non-zero of opposite sign, and 
$\mu_\star(z_0\pm \epsilon)$ is still dominant. 
Then, for all sufficiently large $n$ (depending on $\epsilon$),
$f_n(z_0 \pm \epsilon)$ are real of opposite sign. 
Therefore there exists a root in-between. 

\medskip

In the next section we need some general results on the existence of 
a real sequence of zeros $\{z_n\}$ such that it converges to a real
{\em non-isolated} limiting point. These results can be summarised in the
following  

\begin{lemma} \label{lemma.1}
Assume the hypotheses of Theorem~\ref{theo.BKW}.
Let us suppose that $z_0$ is a real non-isolated limiting point, such
that exactly two dominant eigenvalues $\mu_1$ and $\mu_2$ become 
equimodular at $z_0$. Let us further suppose that:
 \begin{enumerate}
 \item[(a)] The two dominant eigenvalues are analytic functions in a 
            neighbourhood of $z_0$.
 \item[(b)] The eigenvalue $\mu_1$ (resp.\/ $\mu_2$) is dominant (resp.\/ 
       subdominant) in the interval $[z_0-\epsilon,z_0)$, 
       and is subdominant 
       (resp.\/ dominant) in the interval $(z_0,z_0+\epsilon]$, for some 
       $\epsilon>0$. 
 \item[(c)] The corresponding amplitudes $\alpha_1$ and $\alpha_2$ 
       do not vanish at $z_0$.  
 \item[(d)] The eigenvalues $\mu_1,\mu_2$ and the amplitudes 
       $\alpha_1,\alpha_2$ are real in a real neighbourhood of $z_0$. 
 \end{enumerate}
Then:
\begin{enumerate}
 \item If $\alpha_1(z_0) \alpha_2(z_0)>0$ and 
       $\mu_1(z_0)=-\mu_2(z_0)$, then for odd $n$, there is a 
       sequence of real zeros $\{z_n\}$  converging to $z_0$. 
 \item If $\alpha_1(z) \alpha_2(z)<0$ and 
       $\mu_1(z_0)=-\mu_2(z_0)$, then for even $n$, there is a 
       sequence of real zeros $\{z_n\}$ converging to $z_0$. 
 \item If $\alpha_1(z_0) \alpha_2(z_0)<0$ and 
       $\mu_1(z_0)=\mu_2(z_0)$,  then for all $n$, there is a sequence
       of real zeros $\{z_n\}$ converging to $z_0$. 
\end{enumerate}
\end{lemma} 

\bigskip

\noindent
{\bf Remark}. When we say that ``there is a sequence of real zeros 
$\{z_n\}$ converging to $z_0$'', we mean that for all $\epsilon>0$ there 
exists $n_0=n_0(\epsilon)<\infty$ such that for {\em all} $n\ge n_0(\epsilon)$
[or all odd or all even $n$, as the case may be] there is a zero $z_n$ of
$f_n$ satisfying $|z_n-z_0|< \epsilon$.

\medskip 

\proofof{Lemma~\ref{lemma.1}}
The function $f_n$ defined by \reff{def_fn} can be written as
\be
f_n(z) \;=\; \alpha_1(z) \mu_1(z)^n + \alpha_2(z) \mu_2(z)^n +  
          \sum\limits_{j=3}^N \alpha_j(z) \mu_j(z)^n \,, 
\ee
where the $N-2$ other eigenvalues are subdominant in 
a neighbourhood of $z_0$. Then, for every sufficiently small 
$\epsilon>0$ we can find real numbers $0<r_1,r_2<1$ such that 
\be
\left| \frac{\mu_i(z_0-\epsilon)}{\mu_1(z_0-\epsilon)} \right| \;\le\; r_1 \,, 
\ee
for all $i\ge 2$, and 
\be 
\left| \frac{\mu_i(z_0+\epsilon)}{\mu_2(z_0+\epsilon)} \right| \;\le\; r_2 \,, 
\ee
for all $i=1$ or $i\ge 3$. 
We choose $\epsilon$ small enough so that the signs of the dominant 
eigenvalues and amplitudes are the same as at $z_0$. Then, we have that
\begin{subeqnarray}
\slabel{def_fn_minus}
f_n(z_0-\epsilon) &=& \mu_1^n \left[ \alpha_1 + \sum\limits_{j=2}^N 
             \alpha_j \left(\frac{\mu_j}{\mu_1}\right)^n \right] \\
f_n(z_0+\epsilon) &=& \mu_2^n \left[ \alpha_2 + 
             \alpha_1 \left(\frac{\mu_1}{\mu_2}\right)^n + 
             \sum\limits_{j=3}^N 
             \alpha_j \left(\frac{\mu_j}{\mu_2}\right)^n \right] 
\slabel{def_fn_plus}
\end{subeqnarray}
where in \reff{def_fn_minus} [resp.\/ \reff{def_fn_plus}] all quantities are 
evaluated at $z_0-\epsilon$ [resp.\/ $z_0+\epsilon$]. 
Then, for large enough $n$, the quantities $f_n(z_0\pm \epsilon)$  have 
the opposite sign in the following cases:
\begin{enumerate}
 \item $\alpha_1 \alpha_2 > 0$ and $(\mu_1/\mu_2)^n < 0$, which occurs  
  if $\mu_1 \mu_2 < 0$ and $n$ is odd.
 \item $\alpha_1 \alpha_2 < 0$ and $(\mu_1/\mu_2)^n > 0$, which occurs  
  if $\mu_1 \mu_2 > 0$ (then $n$ can have either parity) or if 
  $\mu_1 \mu_2 < 0$ and $n$ is even. 
\end{enumerate}
In any of these cases, the continuous function $f_n$ attains values of 
distinct signs at the endpoints of the interval 
$[z_0-\epsilon,z_0+\epsilon]$,
therefore there should be a zero at some point inside this interval. \qed

\medskip

\noindent
{\bf Remarks.} 1. In the fourth case $\alpha_1\alpha_2 >0$ and 
$\mu_1 = \mu_2$, the zeros converging to $z_0$ are {\em non-real}.

2. If the derivative of the ratio $\mu_1/\mu_2$ is nonvanishing
   at $z_0$, then condition~(b) [or the same condition with
   $\mu_1$ and $\mu_2$ interchanged] necessarily holds.
   However, the converse is false:  it is possible for
   condition~(b) to hold even if $\mu_1/\mu_2$ has a vanishing
   derivative at $z_0$.

\medskip

As remarked in Ref.~\cite[p.~55]{BKW80}, the convergence rate for 
isolated and non-isolated limiting points is rather different: exponentially 
fast for the former $|z_0 -z_n| \le A\, r^n$, and 
$|z_0-z_n| \le A\, n^{-1}$ for the latter, as $n\to\infty$. 
 
%
%
\section{Real zeros of the flow polynomials $\bm{\Phi_{G(nk,k)}}$}
\label{sec.res2}

In this section we will discuss the real zeros of the flow polynomials
$\Phi_{G(nk,k)}$ for $1\le k \le 7$. In particular, we will focus
on the real zeros around $Q=5$, and on the existence of real zeros
$Q>5$.  

%
%
\subsection{$\bm{k\le 5}$}

For $1\le k\le 3$, we compute the flow polynomials of $G(nk,k)$ and their
roots for all $n$ in the range $1\le n\le 30$. 
All the flow roots we have found are smaller than $Q=4$. We conjecture 
that this holds for larger $n$ as well. For $k=4,5$ we find flow polynomials
with real roots greater than $Q=4$:
\begin{itemize}
  \item $G(28,4)$ has two real roots greater than $Q=4$: 
      $Q_1\approx 4.0002086861$ and $Q_2\approx 4.3876416603$.
      As $n$ grows, the maximal real root of $\Phi_{G(4n,4)}$ tends to 
      the value $Q_c(4)\approx 4.5697435537$. 
      The convergence to $Q_0=Q_c(4)$ is due to case (3) of 
      Lemma~\ref{lemma.1}: at this point the dominant eigenvalues come 
      from the sectors with $\ell=1$ and $\ell=3$ links, and  
      $\alpha_{1,(1)}(Q_0) \alpha_{3,(3)}(Q_0) < 0$, and  
      $\mu_{5,1,(1)}(Q_0) = \mu_{5,3,(3)}(Q_0) \approx 11.9477$.

  \item $G(30,5)$ has two real roots greater than $Q=4$: 
      $Q_1\approx 4.0000786673$ and $Q_2\approx 4.4867394006$.
      As $n$ grows, the maximal real root of $\Phi_{G(5n,5)}$ tends to 
      the value $Q_c(5)\approx 4.9029018077$.  
      The convergence to $Q_0=Q_c(5)$ is due to case (2) of
      Lemma~\ref{lemma.1}: at this point the dominant eigenvalues come 
      from the sectors with $\ell=0$ and $\ell=3$ links, and
      $\alpha_{0}(Q_0) \alpha_{3,(3)}(Q_0) < 0$, and  
      $\mu_{5,0}(Q_0) = -\mu_{5,3,(3)}(Q_0) \approx -453.306$.
      Therefore, there are real zeros close to $Q_c(5)$ only for even 
      values of $n$.
\end{itemize} 
Both families thus provide counter-examples to Welsh's conjecture
(Conjecture~\ref{conj.welsh}).  

\bigskip

\noindent
{\bf Remark.}
The values of $Q_c(k)$ mentioned above are obtained by determining the
eigenvalue crossing that corresponds to case~(b) of Theorem \ref{theo.BKW}.

%
%
\subsection{$\bm{k=6,7}$: The isolated limiting point $\bm{Q=5}$}

The family $G(6n,6)$ provides a very strong counter-example to the
Welsh conjecture (Conjecture~\ref{conj.welsh}), as it displays a real 
zero converging to $Q=5$ from below; for instance, $G(144,6)$ has a real 
zero at $Q\approx 4.9987003379$. 

This empirical observation will be made rigorous by applying  
Corollary~\ref{coro.BKW} to the family $G(6n,6)$ at $Q=5$. 
In this case,
there is a unique (and therefore, real) leading eigenvalue, which corresponds
to $\ell=3$ and $\lambda=(3)$: $\mu_{3,(3),\star}(5) \approx 177.122$. 
The corresponding amplitude $\alpha_{3,(3)}$ is given by (\ref{eigen_amp}),
\be
\alpha_{3,(3)}(Q) \;=\; \frac{1}{8} Q(Q-1)(Q-5) \,,
\ee
and is a polynomial in $Q$ with a single zero at $Q=5$. Therefore, the BKW
theorem implies that $Q=5$ is an isolated limiting point for this family.
The fact that the limiting point $Q=5$ is isolated implies that there is an
interval of radius $\epsilon$ around $Q=5$ where there are no other limiting 
points. Therefore, in this interval the eigenvalue $\mu_{3,(3),\star}$ 
is still dominant, and hence real. 
Therefore, Corollary~\ref{coro.BKW} implies that there is a sequence
of {\em real}\/  zeros $\{Q_n\}$ converging to $Q=5$. 
We can go a little bit further and show that there exists a sequence 
of real zeros converging to $Q=5$ {\em from below}:  

\begin{corollary} \label{coro.1}
The point $Q=5$ is an isolated limiting point for the family $G(6n,6)$. 
There is a sequence of real zeros $\{Q_n\}$ of the flow
polynomial $\Phi_{G(6n,6)}$ that converges to $Q=5$ from below. 
\end{corollary}

\proof
Let us consider the point $Q=5-\epsilon$ with $\epsilon$ 
small enough so that: 1) the leading eigenvalue $\mu_{3,(3),\star}(Q)>0$; 
2) all other $N$ sub-leading eigenvalues $\mu_j$ satisfy 
$|\mu_j(Q)/\mu_{3,(3),\star}(Q)|\le r < 1$; 
3) the sub-leading amplitudes $|\alpha_j|\le M$ are bounded; and 
4) $\alpha_{3,(3)}(Q)< 0$. We can always choose $\epsilon$ so that 
these conditions are fulfilled,
as the eigenvalues and amplitudes are analytic functions of $Q$ in a 
neighbourhood of $Q=5$. In addition, we know that $\Phi_{G(6n,6)}(5)>0$,
due to Lemma~\ref{lemma.5flow}.
Therefore, we only need to show that $\Phi_{G(6n,6)}(Q)<0$ for $n$ large
enough. This is easy as, 
\be
\Phi_{G(6n,6)}(Q) \;=\; \mu_{3,(3),\star}(Q)^n \left[ 
    \alpha_{3,(3)}(Q) + 
    \sum\limits_{j=1}^N \alpha_j(Q) 
    \left(\frac{\mu_j(Q)}{\mu_{3,(3),\star}(Q)}\right)^n \right]\,.
\ee
We can always choose $N_0$ such that $NMr^{N_0}<|\alpha_{3,(3)}(5-\epsilon)|$.
Then, for all $n\ge N_0$, the sign of $\Phi_{G(6n,6)}(Q)$ is that of 
$\alpha_{3,(3)}(Q)$ (i.e., negative), so that there should be a point $Q_0$ in
$(5-\epsilon,5)$ such that $\Phi_{G(6n,6)}(Q_0)=0$. \qed 

The same situation applies to the family $G(7n,7)$ at $Q=5$. In this case,
the unique (and real) leading eigenvalue is 
$\mu_{3,(3),\star}(5) \approx -621.779$. 
Therefore, $Q=5$ is an isolated limiting point for this family,
and there is a {\em real}\/ sequence of zeros 
$\{Q_n\}$ converging to $Q=5$. Because $\mu_{3,(3),\star}(5)<0$, the 
above arguments imply that the convergence of the sequence $\{Q_n\}$ 
to $Q=5$ is a bit more complicated:  

\begin{corollary} \label{coro.2}
The point $Q=5$ is an isolated limiting point for the family $G(7n,7)$ 
with $n\ge 3$. There is a sequence of real zeros $\{Q_n\}$ of the flow
polynomial $\Phi_{G(7n,7)}$ converging to $Q=5$. The sub-sequence with
odd $n$ (resp.\/ even $n$) converges to $Q=5$ from above (resp.\/ below). 
\end{corollary}

\proof
If we consider even $n$, then the proof is as before; therefore the
sequence $\{Q_{2p}\}_{p\in\N}$ converges to $Q=5$ from below. 
On the contrary, for odd $n$, then the sign of the leading term 
$\alpha_{3,(3)} \mu_{3,(3),\star}(Q)^n$ is positive for $Q<5$, and
negative for $Q>5$. Therefore, a trivial modification of the above arguments
leads to the convergence of the sub-sequence $\{Q_{2p+1}\}_{p\in\N}$ 
to $Q=5$ from above. \qed 

Corollaries~\ref{coro.1} and~\ref{coro.2} imply 
Theorem~\ref{theo.1}. Moreover, Corollary~\ref{coro.2} implies the existence
of real flow roots arbitrarily close to $Q=5$ on {\em both} sides 
(above and below). Thus, the Haggard--Pearce--Royle conjecture 
(Conjecture~\ref{conj.HPR}) is false, and a counterexample is given by 
$G(7n,7)$ for all sufficiently large odd $n$.

%
%
\subsection{$\bm{k=7}$: Real flow polynomial zeros larger than  $\bm{Q=5}$}
\label{sec.real.flow.roots}
 
There are also real non-isolated limiting points for the families $G(nk,k)$ 
with $1\le k \le 7$. These non-isolated limiting points correspond to
``crossings'' between two dominant eigenvalues, each of them coming from
a block of the transfer matrix $\widehat{\sf T}_{k+1,\ell,\lambda}$ 
with different values of $\ell$. Therefore, both dominant eigenvalues 
are real and analytic in some real interval around the non--isolated
limiting point $Q_c(k)$. In some cases, there is a corresponding sequence 
of real zeros $\{Q_n\}$ converging to that limiting point.

The family $G(6n,6)$ would be in principle a good candidate for having 
real roots larger than $Q=5$: a direct calculation shows that 
$Q_c(6)\approx 5.1079785012$ is a non-isolated limiting point for this family. 
However, the actual computation of all the members of this family up to 
$G(144,6)$ does not reveal any zero $Q>5$. 
The explanation is simple: at $Q=Q_c(6)$ both dominant eigenvalues
are equal 
$\lambda_{1,(1),\star}(Q_c(6))$ $= 
\lambda_{3,(3),\star}(Q_c(6)) \approx 169.757$, and the corresponding 
amplitudes are both positive. Therefore, the hypotheses of 
Lemma~\ref{lemma.1} are not satisfied, and we cannot find {\em real} zeros
converging to the non-isolated limiting point $Q_c(6)$. 
Rather, the zeros converging to $Q_c(6)$ should be non-real.

However, we have found that the family $G(7n,7)$ with {\em odd}\/
$n$ does have members with the desired property. 
Even though $G(63,7)$ does not have any zeros larger than $Q=5$,  
the next member that we have computed, $G(119,7)$, has two such zeros:
$Q_1\approx 5.0000197675$, and $Q_2\approx 5.1653424423$. The flow 
polynomial for $G(119,7)$ has degree $120$ and can be written
as
\be
\Phi_{G(119,7)}(Q) \;=\; (Q-1)(Q-2)(Q-3) P_{117}(Q)\,,
\label{def.PG119.7}
\ee
where $P_{117}(Q) = Q^{117} - 351Q^{116}+ 61191Q^{115} - 
7064107Q^{114} +\ldots$
is a polynomial in $Q$ of degree $117$. The coefficients of this polynomial
are given in Appendix~\ref{appendix.poly}. 
We can formalise the existence of such real roots greater than $Q=5$ in 
the following way: 

\begin{proposition} \label{prop.k=7}
Let $\Phi_G(Q)$ be the flow polynomial of the generalised Petersen graph
$G=G(119,7)$. Then it has a real zero in the interval 
$(5+10^{-5},5+2\times 10^{-5})$, and another real zero in the interval
$(516534 \times 10^{-5},516535 \times 10^{-5})$.
\end{proposition}

\proof
We can evaluate the polynomial $\Phi_G(Q)$ at the two end-points of the
interval $(5+10^{-5},5+2\times 10^{-5})$ using exact rational arithmetic
and find results of distinct sign:
$\Phi_G(5+10^{-5})         \approx +2.21791\times 10^{42}$, and 
$\Phi_G(5+2\times 10^{-5}) \approx -5.27937\times 10^{40}$. Therefore, the
intermediate value theorem ensures the existence of a zero of $\Phi_G$
in the open interval $(5+10^{-5},5+2\times 10^{-5})$. 

The same procedure can be carried out for the second interval:  
$\Phi_G(516534 \times 10^{-5}) \approx -1.46592\times 10^{42}$, 
and 
$\Phi_G(516535 \times 10^{-5}) \approx +4.53729\times 10^{42}$. 
\qed 

\bigskip

\noindent
{\bf Remark.} As a curiosity, 
$\Phi_{G(119,7)}(5)=4488918995790513676672232799446257724715600$,
and $\Phi_{G(119,7)}(4)=1133172760943853528$. 

\bigskip

In fact, the family $G(7n,7)$ has, for every large enough odd $n$, two 
real zeros larger than $Q=5$: one converging to $Q=5$ from above (as 
proven in Corollary~\ref{coro.2}), 
and the other one converging to the limiting point 
$Q_c(7)\approx 5.2352605291$ from below. The complete statement 
of the second part is given by the following result:

\begin{proposition} \label{prop.Qc7} 
The point $Q_c(7)\in (5235260\times 10^{-6},5235261\times 10^{-6})$ 
is a non-isolated limiting point for the family $G(7n,7)$. There is a sequence 
of real zeros $\{Q_n\}$ of the flow polynomial $\Phi_{G(7n,7)}$ converging 
to $Q_c(7)$ from below for odd $n$. 
\end{proposition}
 
\proof
In the interval $[5235260\times 10^{-6},5235261\times 10^{-6}]$,
we only find two dominant eigenvalues
$\lambda_{0,\star}$ and $\lambda_{3,(3),\star}$, each of them coming 
from a different $\ell$ sector. Therefore, both eigenvalues and their
sum $\lambda_{0,\star}+\lambda_{3,(3),\star}$ are  
analytic functions of $Q$ in this interval. 
If we evaluate $\lambda_{0,\star}+\lambda_{3,(3),\star}$ at the two 
end-points of this interval using {\sc Mathematica}, we find 
results of distinct sign: 
$\lambda_{0,\star}(5235260\times 10^{-6})+
\lambda_{3,(3),\star}(5235260\times 10^{-6})\approx 0.000917$, and
$\lambda_{0,\star}(5235261\times 10^{-6})+
\lambda_{3,(3),\star}(5235261\times 10^{-6})\approx -0.000817$. Therefore,
there exists an intermediate value $Q_c(7)$ in the open interval
$(5235260\times 10^{-6},5235261\times 10^{-6})$ such that 
$\lambda_{0,\star}(Q_c(7))+\lambda_{3,(3),\star}(Q_c(7))=0$. Exactly at 
$Q=Q_c(7)$, these two eigenvalues have opposite signs
$\lambda_{0,\star}(Q_c(7)) = 
-\lambda_{3,(3),\star}(Q_c(7)) \approx - 565.833$. 
The corresponding amplitudes are both positive for any $Q>5$. Therefore, 
according to case (1) of Lemma~\ref{lemma.1} we 
find real zeros converging to 
$Q_c(7)$ only for large enough odd values of $n$. 

The fact that the convergence is from below comes from the fact that 
$\alpha_{0,\star}(Q_c(7)) - \alpha_{3,(3),\star}(Q_c(7))>0$, so the
sign of the flow polynomial at $Q_c(7)$ is positive
$\Phi_{G(7n,7)}(Q_c(7))\approx [\alpha_{0,\star}(Q_c(7)) - 
\alpha_{3,(3),\star}(Q_c(7))] \lambda_{0,\star}(Q_c(7))^n > 0$. 
However, its sign at $Q_c(7)-\epsilon$ (for small enough values of $\epsilon>0$)
$\Phi_{G(7n,7)}(Q_c(7)-\epsilon)\approx \alpha_{3,(3),\star}(Q_c(7)) 
\lambda_{3,(3),\star}(Q_c(7))^n < 0$ for odd $n$. Therefore, there should
be a root in-between. \qed

\medskip

Propositions~\ref{prop.k=7} and~\ref{prop.Qc7} imply Theorem~\ref{theo.2}.

%
%
\appendix
\section{Proofs of Lemmas~\ref{lemma.new1}, \ref{lemma.new2}, 
         and~\ref{lemma.degree}} 
\label{appendix.lemmas}

In this appendix we will provide the proofs of the lemmas that are 
essential for proving the main Theorems~\ref{theo.1} and~\ref{theo.2}.

It is useful to rewrite the transfer matrix
${\sf T}_L$ [cf., \reff{TM_decomp}] as the product of two 
operators ${\sf T}_L = {\sf H} \, {\sf V}$, given by  
\begin{subeqnarray}
 {\sf H}   &=& \prod\limits_{i=k}^1 {\sf H}_{0i} {\sf V}_0 \;=\; 
 {\sf H}_{01} {\sf V}_0 {\sf H}_{02} {\sf V}_0 {\sf H}_{03}
               \cdots {\sf V}_0 {\sf H}_{0k} {\sf V}_0 
 \slabel{def_H_full}  \\ 
 {\sf V}   &=& \prod\limits_{i=1}^k {\sf V}_i \;=\; 
               {\sf V}_1 {\sf V}_2 \cdots {\sf V}_{k-1} {\sf V}_k 
 \slabel{def_V_full}
 \label{TM_decomp2}
\end{subeqnarray}
where $L=k+1$, the operators ${\sf V}_i$ and ${\sf H}_{0i}$ are defined
in \reff{def_H_V} with $v=-Q$, and we have used the property 
$[{\sf V}_i,{\sf V}_j]=[{\sf D}_i,{\sf D}_j]=0$ for all $i,j$.  

In this appendix, we will work with a basis consisting on partitions of the 
top-row $\{0,1,\ldots,k\}$, where we have omitted the primes for simplicity. 
An overline over a site means that this site (and the block
it belongs to) is connected to a block of the bottom-row partition by a 
link.  

%
%
\proofof{Lemma~\protect\ref{lemma.new1}} 

Let us consider first the transfer matrix 
$\widetilde{\sf T}_L = {\sf P}_L {\sf T}_L$ with 
${\sf T}_L = {\sf H}\, {\sf V}$ given by 
\reff{TM_decomp2}, and the projector ${\sf P}_L$ defined in \reff{def_PL}. 
Let us fix $\ell$ ($0\le \ell \le k+1$) and a representation 
$\lambda\in S_\ell$ of dimensionality $\dim \lambda$. 
The set $\widetilde{A}_L^{(\ell)}$ can be split into two disjoint sets:  
\begin{itemize}
 \item $\mathcal{C}_1^\ell$ is the set of all partitions of $\{0,1,\ldots,k\}$
       with no un-marked singletons, $\ell$ marked clusters, and vertex $0$
       is a marked singleton.
 \item $\mathcal{C}_2^\ell$ is the set of all partitions of $\{0,1,\ldots,k\}$
       with no un-marked singletons, $\ell$ marked clusters, and vertex $0$
       is not a singleton. In this case, the vertex $0$ can be either 
       marked or un-marked. 
\end{itemize}
It is clear that $|\widetilde{A}_L^{(\ell)}| = 
|\mathcal{C}_1^\ell| + |\mathcal{C}_2^\ell|$. Notice that we 
do not need to consider the case of vertex $0$ being an un-marked singleton,
since by definition partitions of this type do not belong to 
$\widetilde{A}_L^{(\ell)}$. 

In order to find out the structure of the transfer matrix 
$\widetilde{\sf T}_L$, we do not need to specify the representation 
$\lambda$; therefore, we will omit this index (as well as the index $k$) 
to make the notation clearer. 

Let us first consider a partition $\mathcal{P}_1\in \mathcal{C}_1^\ell$
represented by the vector $\bm{e}_{\mathcal{P}_1}$ in the space of partitions. 
Then,
\be
{\sf V} \, \bm{e}_{\mathcal{P}_1} \;=\; (-Q)^k \, \bm{e}_{\mathcal{P}_1} + 
\sum\limits_{\mathcal{P} \in \widetilde{C}_1^\ell} A_{\mathcal{P}} \,  
\bm{e}_{\mathcal{P}} \,,
\label{V_over_C1}
\ee
where the sum is over the set of all partitions $\widetilde{C}_1^\ell$ of 
$\{0,1,\ldots,k\}$ with $\ell$ marked blocks, with vertex $0$ a marked 
singleton, and such that there is at least one un-marked singleton. This
is because ${\sf V}$ only contains the identity operator and detach 
operators ${\sf D}_i$ with $1\le i \le k$. 

If we consider a partition $\mathcal{P}_2\in \mathcal{C}_2^\ell$, then
\be
{\sf V} \, \bm{e}_{\mathcal{P}_2} \;=\; (-Q)^k \, \bm{e}_{\mathcal{P}_2} +
\sum\limits_{\mathcal{P} \in \widetilde{C}_1^\ell} A_{\mathcal{P}} \,
\bm{e}_{\mathcal{P}} + 
\sum\limits_{\mathcal{P} \in \widetilde{C}_2^\ell} B_{\mathcal{P}} \,
\bm{e}_{\mathcal{P}} \;  + 
\sum\limits_{\mathcal{P} \in  C_3^\ell} D_{\mathcal{P}} \, 
\bm{e}_{\mathcal{P}}     \,,
\label{V_over_C2}
\ee
where 
\begin{itemize}

\item $\widetilde{C}_2^\ell$ is the set of all partitions of 
$\{0,1,\ldots,k\}$ with $\ell$ marked blocks, with vertex $0$ not a  
singleton, and such that there is at least one un-marked singleton. 

\item  $C_3^\ell$ is the set of all partitions with $0$ an 
un-marked singleton. But the corresponding terms in (\ref{V_over_C2})
will all be annihilated by the application of the first operator ${\sf V}_0$ 
in ${\sf H}$ \reff{def_H_full}. We can therefore disregard those terms in 
what follows.

\end{itemize}

We now notice that the operator ${\sf H}_{0i}{\sf V}_0$ can be written as:
\be
{\sf H}_{0i}{\sf V}_0 \;=\; 
  -Q \, I +  {\sf D}_0   + Q^2 \, {\sf J}_{0i} 
  -Q \, {\sf J}_{0i} \, {\sf D}_0 
\ee
where ${\sf J}_{0i}$ and ${\sf D}_0$ are the join and detach operators,
respectively. It is not hard to see that for all partitions 
$\mathcal{P}\in \widetilde{A}_{L}^{(\ell)}$:  
\be
{\sf H}_{0i}{\sf V}_0 \bm{e}_{\mathcal{P}} \;=\; -Q \bm{e}_{\mathcal{P}} + 
  \sum\limits_{\mathcal{P}'\in \mathcal{C}_2^\ell} A_{\mathcal{P}'} \,
   \bm{e}_{\mathcal{P}'} \,,
\ee
where if $\mathcal{P}\in \mathcal{C}_2^\ell$, $\mathcal{P}'$ might 
coincide with $\mathcal{P}$. For partitions with at least one un-marked 
singleton $\mathcal{P}\in \widetilde{\mathcal{C}}_1^\ell 
\cup \widetilde{\mathcal{C}}_2^\ell$, we have
\be
{\sf H}_{0i}{\sf V}_0 \bm{e}_{\mathcal{P}} \;=\; -Q \bm{e}_{\mathcal{P}} +
  \sum\limits_{\mathcal{P}'\in \mathcal{C}_2^\ell} A_{\mathcal{P}'} \,
   \bm{e}_{\mathcal{P}'} + 
  \sum\limits_{\mathcal{P}'\in \widetilde{\mathcal{C}}_2^\ell} 
   B_{\mathcal{P}'} \, \bm{e}_{\mathcal{P}'} \,. 
\ee
Notice that when we apply the full operator ${\sf H}$, the partitions
belonging to $\widetilde{\mathcal{C}}_1^\ell 
\cup \widetilde{\mathcal{C}}_2^\ell$ are eliminated by the application of
the operator ${\sf P}_L$ \reff{def_PL}. 

We now put together the above observations, and conclude that if 
$\mathcal{P}\in \mathcal{C}_1^\ell$, then
\be
\widetilde{\sf T}_L \bm{e}_{\mathcal{P}} \;=\; (-Q)^{2k} \bm{e}_{\mathcal{P}} 
 + \sum\limits_{\mathcal{P}'\in \mathcal{C}_2^\ell} A_{\mathcal{P}'} \,
   \bm{e}_{\mathcal{P}'} \,.
\ee
On the other hand, if $\mathcal{P}\in \mathcal{C}_2^\ell$, then
\be
\widetilde{\sf T}_L \bm{e}_{\mathcal{P}} \;=\; 
 \sum\limits_{\mathcal{P}'\in \mathcal{C}_2^\ell} A_{\mathcal{P}'} \,
   \bm{e}_{\mathcal{P}'} \,.
\ee
This means that, if we order the partition states appropriately, the transfer
matrix $\widetilde{\sf T}_L$ has a block triangular form \reff{eq.lemma.new1}.  
This result holds for all representations $\lambda\in S_\ell$. 

The block corresponding to the partitions in $\mathcal{C}_1^\ell$
is diagonal with all eigenvalues equal to $(-Q)^{2k}=Q^{2k}$. 
In terms of the matrix $\widehat{T}_L$ \reff{def_That}, this common eigenvalue 
takes the value $\mu_{k,k+1}= Q^{-2k}\, (-1)^kQ^{2k}= (-1)^k$. \qed 

\medskip

\noindent
{\bf Remarks}. 1. The case $\ell=0$ is special, as $|\mathcal{C}_1^0|=0$; 
 therefore, the diagonal block does not exist.

2. The case $\ell=k+1$ is also special, as $|\mathcal{C}_2^{k+1}|=0$; 
   therefore, the block $\widehat{T}_{k+1,k+1,\lambda}$ is a diagonal matrix 
   with all its elements equal to the trivial eigenvalue $\mu_{k,k+1}$.

\medskip 

%
%
\proofof{Lemma~\ref{lemma.new2}}

Let us first consider the following top-row partition 
$\{\{\overline{0},\overline{1}\},\{\overline{2}\},\ldots,\{\overline{k}\}\}$ 
with $\ell=k$ links. 
This partition does not lead to the trivial eigenvalue $\mu_{k,k+1}=(-1)^k$, 
as it belongs to the class $\mathcal{C}_2^k$ defined in the proof of 
Lemma~\ref{lemma.new1} above. Then the action of 
the transfer matrix ${\sf T}_{k+1}$ [cf., \reff{TM_decomp2}] on this 
partition generates exactly $k$ partitions: 
$\{\{\overline{0},\overline{1}\},\{\overline{2}\},\ldots,\{\overline{k}\}\}$,
$\{\{\overline{1}\},\{\overline{0},\overline{2}\},\ldots,\{\overline{k}\}\}$,
\ldots ,   
$\{\{\overline{1}\},\{\overline{2}\},\ldots,\{\overline{0},\overline{k}\}\}$.
It is important to note that the ordering of the $k$ links is preserved in
the whole process. (Loosely speaking, there is ``no room'' to switch two
links if we have $k+1$ sites and $k$ links.) 
This also means that the action of the transfer matrix ${\sf T}_{k+1}$ 
on any of these partitions is independent of the actual ordering of the links. 
Hence, for each of the $k!$ possible orderings of the $k$ links, the 
transfer matrix will be the same (modulo some reordering of the partitions),
and will have dimension $k$. Therefore, all irreducible representations
of the symmetric group $S_k$ will give the same $k$ non-trivial eigenvalues,
and the multiplicity of each of these eigenvalues in 
$\widehat{\sf T}_{k+1,k,\lambda}$ is just the dimension of the representation
$\lambda$.
\qed 

\medskip

%
%
\proofof{Lemma~\protect\ref{lemma.degree}} 

Let us consider the full matrices $\widehat{\sf{T}}_{k+1,\ell,\lambda}$ and
$\widetilde{\sf{T}}_{k+1,\ell,\lambda}$; 
we will deal with their non-trivial blocks 
$\widehat{\sf{T}}_{k+1,\ell,\lambda}^{(nt)}$ and  
$\widetilde{\sf{T}}_{k+1,\ell,\lambda}^{(nt)}$ at the end. 

All the matrix elements of $\widetilde{\sf{T}}_{k+1,\ell,\lambda}$ are 
polynomials in $Q$. Therefore, 
$\tr(\widetilde{\sf{T}}_{k+1,\ell,\lambda})^n$
will be a polynomial in $Q$. However, not all the matrix
elements of $\widehat{\sf T}_{k+1,\ell,\lambda}$ are polynomials in
$Q$; so we cannot conclude that $\tr (\widehat{\sf{T}}_{k+1,\ell,\lambda})^n$ 
is a polynomial in $Q$.

We need a different argument. Let us consider any partition $\mathcal{P}$
belonging to $\widetilde{\mathcal{A}}_{k+1}^{(\ell)}$. Then, its 
contribution to $\tr (\widetilde{\sf{T}}_{k+1,\ell,\lambda})^n$ will  
be the sum of the contributions of all the diagrams that start with the 
partition $\mathcal{P}$ and end with the same partition after $n$ steps.   
By inspection of the form of the transfer matrix ${\sf T}_{k+1}$ 
\reff{TM_decomp} and its components ${\sf V}_i$ and ${\sf H}_{0i}$ 
\reff{def_H_V}, it is clear that there is a contribution proportional to
$(-Q)^{2kn}$ due to the application of all the identity operators in 
${\sf T}_{k+1}$ (for each layer, there are $3k$ of them, and only $2k$ take 
the factor $-Q$). 
Indeed, this is the minimum number of $-Q$ factors one can possibly 
obtain for a diagram of this type. Therefore, the minimum power of $Q$
that appears in $\tr (\widetilde{\sf{T}}_{k+1,\ell,\lambda})^n$ is $Q^{2kn}$,
which is exactly the inverse of the $n$-th power of the prefactor 
$Q^{-2k}$ in the definition of $\widehat{\sf T}_{k+1}$ 
\reff{Phi_as_trace}/\reff{def_That}. In conclusion, 
$\tr (\widehat{\sf{T}}_{k+1,\ell,\lambda})^n$ is indeed a polynomial in $Q$. 

Once the polynomial character of $\tr (\widehat{\sf{T}}_{k+1,\ell,\lambda})^n$
is established, we have to take care of the degree of this polynomial.
The lemma is proved if we are able to show
that all the entries in $\widetilde{\sf{T}}_{k+1,\ell,\lambda}$ have powers 
of $Q$ of degree at most $3k+\min(1-\ell,0)$. This would imply, that
$\tr (\widetilde{\sf{T}}_{k+1,\ell,\lambda})^n$ is a polynomial in $Q$ of 
degree at most $n[3k+\min(1-\ell,0)]$, and that 
$\tr (\widehat{\sf{T}}_{k+1,\ell,\lambda})^n$ is a polynomial in $Q$ of degree
at most $n[k+\min(1-\ell,0)]$, as claimed. 

Let us first consider the case $\ell=0$. Let us start with an arbitrary
partition of the top row $\mathcal{P}\in \widetilde{\mathcal{A}}^{(0)}_{k+1}$.
When we apply $\widetilde{\sf T}_{k+1}$ to that partition, the maximum number 
of $-Q$ factors that can appear is $(-Q)^{3k}$: we apply the 
part $-Q \bone$ for each vertical operator ${\sf V}_i$, and the join 
operator $-Q {\sf J}_{0i}$ for each horizontal operator ${\sf H}_{0i}$. 
We obtain for all $\mathcal{P}\in \widetilde{\mathcal{A}}^{(0)}_{k+1}$
the same final partition: $\{\{0,1,\ldots,k\}\}$. Therefore, 
$\tr (\widetilde{\sf{T}}_{k+1,0,\lambda})^n$ is a polynomial in $Q$ of degree
at most $3kn$. 

The case $\ell=1$ is similar. If we start with an arbitrary
partition of the top row $\mathcal{P}\in \widetilde{\mathcal{A}}^{(1)}_{k+1}$ 
with $0$ marked, and we apply the same operators as above, we end up with
the partition $\{\{\overline{0},\overline{1},\ldots,\overline{k}\}\}$. 
The coefficient corresponding to this partition is a polynomial 
in $Q$ of degree at most $3k$. 
Therefore, $\tr (\widetilde{\sf{T}}_{k+1,1,\lambda})^n$ is a polynomial 
in $Q$ of degree at most $3kn$. 

The case $2\le \ell \le k+1$ is similar. Let us start with the following
simple partition in $\widetilde{\mathcal{A}}^{(\ell)}_{k+1}$: 
$\{\{\overline{0},\overline{\ell},\ldots,\overline{k}\},
\{\overline{1}\}, \ldots, \{\overline{\ell-1}\}\}$. 
The argument is similar to the previous cases, except that we cannot 
join the site $0$ to any of the other $\ell-1$ blocks; since if we did, 
then the number of links would be smaller than $\ell$. Therefore, we get a
diagonal entry which is a polynomial in $Q$ of degree at most 
$3k-(\ell-1)$. Therefore, $\tr (\widetilde{\sf{T}}_{k+1,\ell,\lambda})^n$ is 
a polynomial in $Q$ of degree at most $n[3k-(\ell-1)]$.

Finally, Eq.~\reff{def_reduced_traces} implies that 
$\tr (\widehat{\sf{T}}_{k+1,\ell,\lambda})^n$ and 
$\tr (\widehat{\sf{T}}_{k+1,\ell,\lambda}^{(nt)})^n$ differ by a term of
order $Q^0$. Therefore, 
$\tr (\widehat{\sf{T}}_{k+1,\ell,\lambda}^{(nt)})^n$ is a polynomial in 
$Q$ of degree at most $n[k+\min(1-\ell,0)]$. 
\qed

\section{More structural properties of the transfer matrix}
\label{sec.extra}

%
%
Lemma~\ref{lemma.new1} proved that the relevant diagonal blocks 
$\widehat{\sf T}_{k+1,\ell,\lambda}$ have an upper-block-triangular
structure \reff{eq.lemma.new1}. One interesting question is to know the
dimension of the diagonal block $\widehat{\sf D}_{k+1,\ell,\lambda}$. 
This is the content of the following lemma:

\begin{lemma}
\label{lemma.mult.trivi}
Fix $k\ge 1$, and $\ell=0,\ldots,k+1$. Then, the dimension of the 
trivial diagonal block $\widehat{\sf D}_{k+1,\ell,\lambda}$ 
can be written as 
\be
\dim \widehat{\sf D}_{k+1,\ell,\lambda} \;=\; 
     \widetilde{N}_{k,1}(\ell) \dim\lambda\,,
\ee
where $\widetilde{N}_{k,1}(\ell)$ is given by 
\be
\widetilde{N}_{k,1}(\ell) \;=\; \begin{cases}
    0 & \quad \text{if $\ell=0$} \\
    \displaystyle
    \frac{1}{(\ell-1)!} \, |\widetilde{A}_k^{(\ell-1)}|
    \;=\; \sum\limits_{p=0}^k \binom{k}{p} \stirlingsubset{p}{\ell-1}
    S_{k-p} & \quad \text{if $1\le \ell \le k+1$}
   \end{cases}
\label{def_Ntildek1}
\ee
and $S_n$ is given in \reff{def_Sn}.
\end{lemma}

\proof
As we have seen in the first remark after the proof of Lemma~\ref{lemma.new1},
$|\mathcal{C}_1^0|=0$, so $\dim \widehat{\sf D}_{k+1,0,\lambda}=0$. 
Thus, let us  assume that $1\le \ell\le k+1$. The dimension of the diagonal 
block $\widehat{\sf D}_{k+1,\ell,\lambda}$
should be equal to the number of partitions of the set 
$\{1,2,\ldots,k\}$ with no un-marked singletons, $\ell-1$ marked blocks,
and corresponding to the representation $\lambda\in S_\ell$ of 
dimensionality $\dim\lambda$.
But now, as we are considering linear combinations of the 
partitions with the right symmetries under $S_\ell$, these blocks should
be considered indistinguishable. Therefore,  
\be
\dim \widehat{\sf D}_{k+1,\ell,\lambda} \;=\; 
\frac{1}{(\ell-1)!} \, |\widetilde{A}_{k}^{(\ell-1)}| \, 
      \dim \lambda \,. 
\ee
The total number of trivial eigenvalues in 
$\widehat{\sf T}_{k+1,\ell,\lambda}$ divided by $\dim\lambda$ is 
\be
\widetilde{N}_{k,1}(\ell) \;=\; 
\frac{1}{(\ell-1)!}  \, |\widetilde{A}_{k}^{(\ell-1)}| \,. 
\ee
Using \reff{def_Atilde_kl} and after some algebra, it is not
difficult to find the expression \reff{def_Ntildek1}. \qed  

\medskip

\noindent
{\bf Remark.} Indeed, the dimensions of the matrices 
$\widehat{\sf D}_{k+1,\ell,\lambda}$ exactly found by a computer-assisted proof
in Section~\ref{sec.res1.gen} do coincide with the analytic formula above.    

\bigskip

The next step is to compute the dimensionality of the non-trivial block 
$\widehat{\sf T}_{k+1,\ell,\lambda}^{(nt)}$. 

%
%
\begin{lemma} \label{lemma.dim.ntK}
Fix $k\ge 1$, $\ell \in \{0,1,\ldots,k+1\}$, and the representation 
$\lambda\in S_\ell$. Then the dimension of the non-trivial diagonal block  
$\widehat{\sf T}_{k+1,\ell,\lambda}^{(nt)}$ 
\be
\dim \widehat{\sf T}_{k+1,\ell,\lambda}^{(nt)} \;=\; 
\widetilde{N}_{k,0}(\ell,\lambda)\,,
\ee
is given by 
\be
\widetilde{N}_{k,0}(\ell,\lambda) \;=\; \begin{cases}
  \displaystyle
  S_{k+1} + \sum\limits_{p=0}^{k} \binom{k+1}{p+1} S_{k-p} & \quad
  \text{if $\ell=0$} \\[2mm]
  \displaystyle
  \dim\lambda \,  \sum\limits_{p=1}^k \binom{k}{p} \stirlingsubset{p}{\ell}
    \left[S_{k-p} + S_{k+1-p}\right] & \quad
  \text{if $1\le \ell \le k-1$} \\[2mm]
   k \dim \lambda     &  \quad\text{if $\ell =k$} \\
   0    &  \quad\text{if $\ell =k+1$}
 \end{cases}
\label{def_Nk0_OK}
\ee
and $S_n$ is given by \reff{def_Sn}. 
\end{lemma}

\proof
Let us first fix $\ell$ ($0\le \ell \le k+1$) and the 
irreducible representation $\lambda \in S_\ell$ of dimensionality $\dim\lambda$.
The dimension of the corresponding diagonal block 
$\widehat{\sf T}_{k+1,\ell,\lambda}$ is 
\be
\dim \widehat{\sf T}_{k+1,\ell,\lambda} \;=\;  
\frac{\widetilde{N}_k(\ell)}{\ell!} \, \dim\lambda \,.
\label{def_dim_T_L_l_lambda}
\ee
The dimension of the trivial block ${\sf D}_{k+1,\ell,\lambda}$ is  
$\widetilde{N}_{k,1}(\ell) \dim \lambda$ by Lemma~\ref{lemma.mult.trivi}. 
Therefore, the number of non-trivial eigenvalues in this block is 
given by 
\be
\dim \widehat{\sf T}_{k+1,\ell,\lambda}^{(nt)} \;=\;  
     \widetilde{N}_{k,0}(\ell,\lambda) \;=\; 
\left[\frac{\widetilde{N}_k(\ell)}{\ell!}  \,-\, 
            \widetilde{N}_{k,1}(\ell) \right]\, \dim\lambda \,.
\label{def_tildeN_k0}
\ee

We can provide a closed form for $\widetilde{N}_{k,0}(\ell,\lambda)$ 
for $0\le \ell\le k-1$ by combining Lemma~\ref{lemma.mult.trivi}, 
the definition $\widetilde{N}_k(\ell)=|\widetilde{A}_{k+1}^{(\ell)}|$
\reff{def_Atilde_kl}, and \reff{def_tildeN_k0}. After some algebra we
find Eq.~\reff{def_Nk0_OK} for $\ell\le k-1$.

The case $\ell=k$ is derived directly from Lemma~\ref{lemma.new2}. 
There are $k$ distinct non-trivial eigenvalues for 
every irreducible representation $\lambda$ of $S_k$, each of them with
multiplicity $\dim\lambda$. 
Therefore, $\widetilde{N}_{k,0}(k,\lambda)=k \, \dim\lambda$. 

Finally, for $\ell=k+1$, the number of non-trivial eigenvalues is zero,
as all eigenvalues are trivial in this sector 
(see the second remark after the proof of Lemma~\ref{lemma.new1}).  
\qed   

Lemmas~\ref{lemma.new1} and~\ref{lemma.new2}
imply that the flow polynomial \reff{def_Phi_G} for the generalised Petersen 
graph $G(nk,k)$ is given by the following ``complete'' decomposition for 
any $k\ge 1$:
\be
\Phi_{G(nk,k)}(Q) \;=\; \sum\limits_{\ell=0}^{k-1}
   \sum\limits_{\lambda\in S_\ell} \alpha_{\ell,\lambda}
   \sum\limits_{s=1}^{\widetilde{N}_{k,0}(\ell,\lambda)} 
                                            \mu_{k,\ell,\lambda,s}^n 
   + \beta_k \sum\limits_{s=1}^k \mu_{k,k,s}^n
   + \gamma_{k+1} (-1)^{nk} 
\label{eq.phi}
\ee
where $\beta_\ell$ is given in \reff{eigen_amp_total}, and  
$\gamma_{k+1}$ is given by 
\be
\gamma_{k+1} \;=\; \beta_{k+1} + 
                 \sum\limits_{\ell=1}^k \sum_{\lambda \in S_\ell} 
                 \alpha_{\ell,\lambda} \widetilde{N}_{k,1}(\ell) \dim\lambda 
             \;=\; \beta_{k+1} + 
                 \sum\limits_{\ell=1}^k \beta_\ell \widetilde{N}_{k,1}(\ell)\,. 
\label{def.gamma}
\ee

\medskip

\noindent
{\bf Remark.} Lemma~\ref{lemma.new4} implies that all eigenvalues appearing 
in \reff{eq.phi} are distinct for $1\le k\le 7$. We conjecture that
all eigenvalues in \reff{eq.phi} are also distinct for each $k\ge 8$.
 
\medskip

We can compute the coefficients $\gamma_{k+1}$ by using \reff{def.gamma}, the 
expressions \reff{def.values.beta} for the amplitudes $\beta_\ell$, and the
values \reff{def_Ntildek1} of $\widetilde{N}_{k,1}(\ell)$. 
The results for $1\le k\le 7$ are:
\begin{subeqnarray}
\gamma_2 &=& Q^2 - 3Q^2 +   1                             \\
\gamma_3 &=& Q^3 - 5Q^2 +  6Q   -   1                     \\
\gamma_4 &=& Q^4 - 7Q^3 + 15Q^2 -  11Q   +   1            \\
\gamma_5 &=& Q^5 - 9Q^4 + 28Q^3 -  38Q^2 +  20Q   -  1    \\
\gamma_6 &=& Q^6 -11Q^5 + 45Q^4 -  90Q^3 +  90Q^2 - 27Q   +   1          \\
\gamma_7 &=& Q^7 -13Q^6 + 66Q^5 - 175Q^4 + 260Q^3 -207Q^2 +  70Q   - 1   \\
\gamma_8 &=& Q^8 -15Q^7 + 91Q^6 - 301Q^5 + 595Q^4 -707Q^3 + 469Q^2 -135Q +1
\label{def.values.gamma}
\end{subeqnarray}
For instance, for $k=3$, we have that 
$\gamma_4=\beta_1+4\beta_2+3\beta_3+\beta_4$.

For $1\le k\le 7$, the total number of {\em distinct} eigenvalues is given by
\be
\widetilde{D}_k \;=\; 1 + k + \sum\limits_{\ell=0}^{k-1} 
                              \sum\limits_{\lambda\in S_\ell} 
                              \widetilde{N}_{k,0}(\ell,\lambda)\,,
\ee
where the $\widetilde{N}_{k,0}(\ell,\lambda)$ are given in 
\reff{def_Nk0_OK}. The
values are $\widetilde{D}_k=3,7,36,229,1658,12803,105934$ for $k=1,\ldots,7$. 

\clearpage
%
%
\section{The polynomial $\bm{\Phi_{G(119,7)}}$} \label{appendix.poly}
%
%
\begin{table}[h]
\centering
\vspace*{-1cm}
\begin{tabular}{rr}
\hline\hline
\multicolumn{1}{c}{$i$} & \multicolumn{1}{c}{$a_i$} \\ 
\hline
0 &240453758183717079931230416441214627161100181583221695778758847017660\\
1 &3778010581676303383947166862404894626185168386864045862610150115812052\\
2 &30335569899732630785396756910315613411372095140384361997526113182374163\\
3 &165825101972051263346220423069919998872239305113596582973635632055394977\\
4 &693579540994844285783536577772823633624035569377766733899721814879025261\\
5 &2365298420910039361031441041999881902182643546761840393901756834507719194\\
6 &6843881667041711337123509683076504016112680686224789087846006938716709936\\
7 &17263457862032800683223779137831846913199242114984324936441639648268053289\\
8 &38713433713214705146001095156944656719927970414458008688062316276391488388\\
9 &78323139898934533680532026518046626525971702819240021820160355852376101900\\
10 &144596924706446634397013884950643553134413962554815223156161387400333508225\\
11 &245805490880525223783583677086377232647690995050227537846281841251646914824\\
12 &387578029394911704914155641300482722338098395075782035706404804765273096365\\
13 &570251442350241862197997781738305884386349966454124613787648726411773264858\\
14 &786830079084979854704322218255681568053263515605219068314256695784834361687\\
15 &1022415665287852311682199189937274383849816815595817077282491287814210089327\\
16 &1255617405561053292049115611630960519206807161888631329237623843787650902598\\
17 &1461844034312439449575494091690068905133069052914279852509149780627087441657\\
18 &1617740821249798765950531657277432254427341252461560806920840933023592630423\\
19 &1705621288663026987807981694789591864355698512752538498113371227662898599836\\
20 &1716721077841948079359918580031477807370750494656542267873960267539031216693\\
21 &1652471749735377739204603329938965183162646296909313884512941374410264086521\\
22 &1523600861689858982867164652432211354402872413669859092087649541160917763276\\
23 &1347478648579683679991904590350043152972817068987716895096810710111624258414\\
24 &1144543577841773436389066666439373326178625474359360165078175860151437888418\\
25 &934744903425435172192816740334988586706047594412696933301815485640045053088\\
26 &734761463017109267283189805699271455927373168195072916098687225728186938239\\
27 &556406139434285183296690956661070857540502747488343675526234531927677714022\\
28 &406248919977322871211475078592598791327840251152453917692944363153498710100\\
29 &286204685555683178005609409552066982429730669340777755008844846916076522994\\
30 &194690835588052864277724799113473665492499996408415866330514336755095511519\\
\hline\hline
\end{tabular}
\caption{
Coefficients of the flow polynomial
$\Phi_{G(119,7)}(Q) = (Q-1)(Q-2)(Q-3) \sum\limits_{i=0}^{117} 
(-1)^{i+1} a_i Q^i$.
}
\label{table.P117}
\end{table}

\clearpage
\thispagestyle{empty}
\addtocounter{table}{-1}
\begin{table}[htb]
\centering
\begin{tabular}{rr}
\hline\hline
\multicolumn{1}{c}{$i$} & \multicolumn{1}{c}{$a_i$} \\
\hline
31 &127959409173717651616971568749651845945079328780268936837379415435944126802\\
32 &81303058361075417592304554906219204545070829443708169095378279034181616797\\
33 &49966395799032436266189002396534646533235474102922584729978792549768709570\\
34 &29716297381875173070101426039892101696811371497930415707024934873184994971\\
35 &17109926011423464771467167427507814898291755276949983853151230858292134005\\
36 &9541424411645170221820026129682497673840857753028326218052776323947917867\\
37 &5155250259203698010398273841766642027281260863127250587163264787844562005\\
38 &2699625968148386503782707617994694138490854643501950268752229332983496484\\
39 &1370590172276847791788823798202187745155060493411059843973697377813429115\\
40 &674811604481186691941181395279228274700479647745836519066846509110686871\\
41 &322283263195622661124593258235203925484545860053525952388302407768566612\\
42 &149338540704756141327132733140880585247255726988170348306522536601036428\\
43 &67154254327764896250779574488797365684816293909986077788680077649219059\\
44 &29310387712861712336349988437753297027267572645934195793910184667409525\\
45 &12418975519744799441679142496250486367091770751393334192303593690937216\\
46 &5108897246470415157235242805025854028977310258147102091244333174463070\\
47 &2040790191809533248688898356534425242558953177861065156649563146370679\\
48 &791669367678994305215119572575588439395174310193863932545018547493949\\
49 &298263611514781343272224652335166519468558609433860386772388696144855\\
50 &109143202315633324987330964029233262735802381917881292972386146290583\\
51 &38793120529207946517836936596295050010421179961661603601012761033976\\
52 &13393333822523527127943242411311403493467668611668190366289160119798\\
53 &4491649304261905729281333948467503543032844169847443690259485449787\\
54 &1463208911995358578670442816667443146618570037914577188051167519252\\
55 &463002576567921398643877003491825254798116352410207026582354227639\\
56 &142306049793794840801481088729855032253716088115496415256458286993\\
57 &42482066138212698853729014913405235508231962094892775535623301954\\
58 &12316937820630034457717661311813898450293797908401593794293724758\\
59 &3468010035454592592294532401090077198723295654493559716526595085\\
60 &948193746090263874793189178585310649056926723273402296835255341\\
61 &251712395999861219430022121008687006528308619159639108633408792\\
62 &64870522436154805702181162937007326080100049608985554185883024\\
63 &16227914273399692643538894326321813038077624088578343763024693\\
64 &3939856175099340812778271579647788724054817562072496203271622\\
65 &928159554704604760538186427442530041159446537967362956747430\\
66 &212130925228863323822024721431112601866613699451131860235710\\
67 &47025221622368034908085340305586428332317539464533841648123\\
68 &10108841512034189983877166351043431883332647922883010856327\\
69 &2106711947665769238334508701464763203793865636828903442294\\
70 &425522761043753110078251752857561183613194192306420512412\\
71 &83277105914548332097899604836384635641176843712157856926\\
72 &15786114232928385188569755800951974220820338221822797673\\
73 &2897497349085917262770771233534175446735608061805601330\\
74 &514765848404059046836837761417162573532059958487097353\\
75 &88483984104925516580422896030293322509203586666446608\\
\hline\hline
\end{tabular}
\caption{
(Continued.)
}
\label{table.P117B}
\end{table}

\clearpage
\addtocounter{table}{-1}
\begin{table}[htb]
\centering
\begin{tabular}{rr}
\hline\hline
\multicolumn{1}{c}{$i$} & \multicolumn{1}{c}{$a_i$} \\
\hline
76 &14709751385762293879136933265239581784434168823726105\\
77 &2363940443321701747333111550458096052434867769928216\\
78 &367070636910713149669510355126240632283426228583990\\
79 &55045491030169621973045797552387375721405805002210\\
80 &7967363593383576742789482483060195429602666353792\\
81 &1112432655029526956212696750602542449173138786358\\
82 &149736221917777447936668185649515326795606657139\\
83 &19417037687233804941366229781880098811952696949\\
84 &2423988165427613282917222688678473423362471885\\
85 &291095664115059080227254369870423630080267269\\
86 &33600127040413356951575245326064861734434090\\
87 &3724448845910258224439985334986604271814630\\
88 &396083067053161605214005007534380917557302\\
89 &40370994515077502586895936288134103918883\\
90 &3939422436008602219258246065351204112198\\
91 &367585380265857178359590490307505526698\\
92 &32755670407906919539877604731444999292\\
93 &2783619278733903044429223967432069173\\
94 &225251413853620678004411181186579261\\
95 &17327620037637857872926704478452708\\
96 &1264839938089389154908070831758552\\
97 &87436328017906227077510426235281\\
98 &5711550730808768130994711119064\\
99 &351695047292333601570905687819\\
100 &20358915651379087237820736779\\
101 &1104604810216428455300363966\\
102 &55981869638360027493471673\\
103 &2640008664516116400095478\\
104 &115338644735668374521021\\
105 &4644751096461226418781\\
106 &171402937526338140238\\
107 &5756292666737489361\\
108 &174485239351077202\\
109 &4726362454421551\\
110 &113000004047458\\
111 &2347463061581\\
112 &41510490693\\
113 &607499109\\
114 &7064107\\
115 &61191\\
116 &351\\
117 &1\\
\hline\hline
\end{tabular}
\caption{\
(Continued.)
}
\label{table.P117C}
\end{table}

\clearpage
\section*{Acknowledgments}

We wish to thank Gordon Royle for suggesting us to study the generalised
Petersen graphs, and Alan Sokal for helpful conversations and for a 
critical reading of the manuscript that gave rise to many constructive
comments. We also wish to thank an anonymous referee for his/her 
very detailed suggestions and criticisms that helped us to improve 
the manuscript.

Both authors also thank the Isaac Newton Institute for Mathematical
Sciences, University of Cambridge, for hospitality during the
programme on Combinatorics and Statistical Mechanics (January--June
2008), where this project started. J.S. also thanks the \'Ecole Normale
Sup\'erieure for hospitality in June 2009, June 2011, and December
2012.

The research of J.S. was supported in part by Spanish MICINN/MINECO grants
MTM2008-03020, FPA2009-08785, MTM2011-24097 and FIS2012-34379,
and by U.S.\ National Science Foundation grant PHY--0424082. 
The research of J.L.J. was supported in
part by the European Community Network ENRAGE (grant
MRTN-CT-2004-005616), and by the Agence Nationale de la Recherche
(grant ANR-06-BLAN-0124-03).

\end{document}